\documentclass[review,onefignum,onetabnum]{siamart190516}
 
\usepackage{lipsum}
\usepackage{amsfonts}
\usepackage{epstopdf}
\usepackage{algorithmic}
 
\usepackage{float}
\usepackage{graphicx}

\usepackage{caption,subcaption}

\ifpdf
  \DeclareGraphicsExtensions{.eps,.pdf,.png,.jpg}
\else
  \DeclareGraphicsExtensions{.eps}
\fi

\graphicspath{ {./images/} }

\newsiamremark{remark}{Remark}
\newsiamremark{hypothesis}{Hypothesis}
\crefname{hypothesis}{Hypothesis}{Hypotheses}
\newsiamthm{claim}{Claim}
\newtheorem{example}{Example}[section]
\headers{NTHP for Sparse LCP via A New Merit Function}{S.L. Zhou, M.J. Shang, L.L. Pan and M. Li }
 
\title{Newton Hard Thresholding Pursuit for Sparse LCP via A New Merit Function\thanks{Submitted to the editors DATE.
\funding{This work is supported in part by the National
Natural Science Foundation of China (11601348, 11801325, 11771255, 11971052), ``111" Project of China (B16002) and Young Innovation Teams of Shandong Province (2019KJ1013).}}}

\author{Shenglong Zhou\thanks{School of Mathematics, University of Southampton, Southampton SO171BJ, United Kingdom 
  (\email{shenglong.zhou@soton.ac.uk}).}
\and  Meijuan Shang\thanks{Corresponding author. College of Science, Shijiazhuang University, Shijiazhuang 050035, People's Republic of China (\email{shangmj1108@163.com}).}
\and Lili Pan \thanks{Department of Mathematics, Shandong University of Technology, Zibo 255000, People's Republic of China  (\email{panlili1979@163.com}).  }
   \and Mu Li\thanks{Mechatronics Engineering Group, University of Southampton, Southampton SO171BJ, United Kingdom (\email{m.li@soton.ac.uk}).}
}
\usepackage{amsopn}

\def\R{{\mathbb  R}}
\def\N{{\mathbb N}}
\def\T{{\cal T}}

\def\bx{ {x}}

\def\sol{{\texttt{sol}(M,q)}}

\def\fea{{\texttt{fea}(M,q)}}
\def\feas{{\texttt{fea}$_s(M,q)$}}

\def\NHTP{\texttt{NHTP}}
\def\HNHTP{\texttt{HNHTP}}
\def\HTP{\texttt{HTP}}
\def\ETA{\texttt{ETA}}

\def\Diag{{\rm Diag}}
\def\supp{{\rm supp}}

\ifpdf
\hypersetup{
  pdftitle={NHTP for Sparse LCP via A New Merit Function},
  pdfauthor={S.L. Zhou,  M.J. Shang, L.L. Pan and M. Li  }
}
\fi

\nolinenumbers
\begin{document}

\maketitle

\begin{abstract}
{{Solutions to the linear complementarity problem (LCP) are naturally sparse in many applications such as bimatrix games and portfolio section problems. Despite that it gives rise to the hardness, sparsity makes optimization faster and enables relatively large scale computation.  Motivated by this, we take the sparse LCP into consideration, investigating the existence and boundedness of its solution set as well as introducing a new merit function, which allows us to convert the problem into a sparsity constrained optimization.}} The  function turns out to be continuously differentiable and twice continuously differentiable for some chosen parameters. Interestingly, it is also convex if {{the involved matrix}} is positive semidefinite. We then explore the relationship between the solution set to the sparse LCP and stationary points of the sparsity constrained optimization. Finally, Newton hard thresholding pursuit is adopted to solve the sparsity constrained model. Numerical experiments demonstrate that the problem can be efficiently solved through the new merit function.
\end{abstract}


\begin{AMS}
   {{	90C33, 90C2, 90C30}}
\end{AMS}

\section{Introduction}

The linear complementarity problem (LCP) aims at finding a vector $x\in\R^n$ such that
\begin{eqnarray}\label{CP}
x\in\sol:= \{x\in\R^n:~x\geq0,\  Mx+q \geq0, \ \langle x, Mx+q\rangle=0 \},
\end{eqnarray}
where $M\in\R^{n\times n}$ and $q\in\R^{n}$. Here, $x\geq0$ {{means that each element}} of $x$ is nonnegative. Linear complementarity problems have extensive applications in economics and engineering {{such as Nash equilibrium problems, traffic equilibrium problems, contact mechanics problems and option pricing, to name a few. More applications can be found in \cite{CPS92,FP03,FMP01} and the references therein.  Among them, there is an important class  trying to seek for a solution where most of its elements are zeros, namely, a sparse solution. For example,  players in bimatrix games are willing to choose a small portion of reasonable strategies from a set of pure strategies to save their computational time. In the portfolio selection problem, most investors are only interested in a `small' portfolio  from a group of assets, see more details in \cite{CPS92,XHZ08, SZX14}. Mathematically, these examples can be characterized as the following sparse LCP}}
\begin{eqnarray}\label{sncp}
x\in\sol \cap S\ \ \ \ {\rm with}\ \ \ \ S:=\{x\in\R^n:~\|x\|_0\leq s\},
\end{eqnarray}
where $\|x\|_0$ is the zero norm of $x$, which counts the number of nonzero elements of $x$, and $s\ll n$ is a positive integer. Note that $\|\cdot\|_0$ is {{not a norm}} in the sense of the standard definition.  In order to address the LCP, a commonly used approach is {{to convert the problem into an unconstrained minimization problem through the  NCP (nonlinear complementarity problem) functions. A function $\psi:\R\times\R\rightarrow \R$ is called an NCP function  if it satisfies}}
\begin{eqnarray} \psi (a,b)=0\ \  \Longleftrightarrow\ \ a\geq0,\ \ b\geq0, \ \ ab=0.\end{eqnarray}
In this paper, we introduce a new function {{$\phi_r:\R\times\R\rightarrow \R$}}  defined by
\begin{eqnarray}\label{func-phi} \phi_r(a,b):=\frac{1}{r}\Big[a_+^rb_+^r+(-a)_+^r+(-b)_+^r\Big]:=\frac{1}{r}\Big[a_+^rb_+^r+|a_-|^r+|b_-|^r\Big].\end{eqnarray}
where $r>0$ is a given parameter, $a_+:=\max\{a,0\}$ and $a_-:=\min\{a,0\}$.
It is easy to see that $\phi_r$ is indeed an NCP function for any given $r>0$. However, through this paper, we only focus on choices of $r\geq2$. Because this new  function is proven to be continuously differentiable everywhere for any $r\geq2$ and  twice continuously differentiable for any $r>2$, see  \Cref{pro-phi-1}. When it comes to model (\ref{CP}), we construct a new merit function $f_r$ through $\phi_r$  as
\begin{eqnarray}\label{sco-obj}
f_r(x)&:=&\sum_{i=1}^n \phi_r(x_i,M_ix+q_i)\\
&=&\frac{1}{r}\Big[ \left\langle x_+^r, (Mx+q)_+^r\right\rangle +\|x_-\|_r^r+\|(Mx+q)_-\|_r^r\Big],\nonumber
\end{eqnarray}
where $\|x\|_r^r:=\sum_i |x_i|^r$ (particularly, write $\|\cdot\|:=\|\cdot\|_2$), $M_i$ is the $i$th row of $M$ and $x_+^r$ and  $x_-$ are defined by (\ref{notation}). Clearly,   $f_r(x)\geq0$ for any $x\in\R^n$. Based on this function, {{to solve the sparse LCP (\ref{sncp}) for a given $s\in\N$ and $s\ll n$}},  we will   deal with the following sparsity constrained optimization throughout this paper
\begin{eqnarray}\label{slco}
\min_x  \ f_r(x),\ \ \ \ {\rm s.t.}\ \  x\in S.
\end{eqnarray}

\subsection{NCP functions} There are numerous NCP functions that have been proposed. One of the most well-known functions is the Fischer-Burmeister (FB) function. It was first introduced by Fischer in \cite{F92} and widely used in designing semismooth Newton type methods for solving mathematical programming with   complementarity conditions. {{Then many  variants have been investigated, see \cite{KK98,LW09} and \cite{CP08} for more information. All those functions share a similar mathematical formula and hence enjoy similar properties. They are continuously differentiable everywhere except at the origin where their Hessians are unbounded.  In \cite{CH97}, the authors}} took advantage of the natural residual (namely, minimum function)  to construct an NCP function, with a simple structure but offering little of the second order information. It is continuously differentiable everywhere as well but nondifferentiable at the origin and along a line. {{The authors in \cite{CCK00} cast an NCP function through the convex combination of the FB-function and the maximum function. The function is continuously differentiable everywhere except at the solution set \eqref{CP}. In \cite{MS93}, a continuously differentiable implicit Lagrangian, an NCP function, was explored.}} Another interesting class of functions have been studied by authors in \cite{KYF97}. They are able to be twice continuously differentiable if their involved parameters are chosen properly.   Functions mentioned above have drawn much attention and have been shown to enjoy many favourable properties  \cite{MS93,F95,YF95,F96, GK96,K96,M96,KYF97,LT97,steffensen2010new,pan2014generalized}.

\subsection{Contributions} Contributions of this paper are summarized below.
\begin{itemize}
\item[i)] We propose a new type of NCP function $\phi_r$, which allows us to  construct a new merit function $f_r$ to deal with the LCP.  It turns out that $f_r$ is continuously differentiable everywhere for any $r\geq2$ and  twice continuously differentiable for any $r>2$, see  \Cref{pro-LCP}. Moreover, if the matrix $M$  is positive semidefinite, then $f_r$ is convex. This means, in order to solve the LCP,  one could address an unconstrained convex optimization that minimizes $f_r$, {{namely, find a stationary point of $f_r$ which by the convexity is a solution to ${\rm min}_x f_r(x)$}}. We then reveal the relationship  between a solution to the LCP and a stationary point, see  \Cref{the-sta}. 

\item[ii)] Not only do we prove the existence and the boundedness of the solution set to the sparse LCP, and the boundedness of the level set of  $f_r$ over $S$, but we also establish the relationship between a solution  to the sparse LCP and a stationary point to the sparsity constrained optimization (\ref{slco}).

 \item[iii)] To process the sparsity constrained optimization  (\ref{slco}), we take advantage of the Newton hard thresholding pursuit (\NHTP) method proposed in \cite{ZXQ19}, whose convergence results are well established in \Cref{sec:nhtp}. Numerical experiments demonstrate that the adopted  method has excellent performance to solve the sparse LCP in terms of the fast computational speed and high order of accuracy. {{What is more, we apply the method  to deal with (\ref{slco}), where the merit objectives are constructed from three existing famous NCP functions.}} Numerical comparisons show that \NHTP\ performs much better on solving the model with $f_r$  than solving models with the other merit functions. {{In a nutshell, the sparse LCP can be solved more effectively by converting it into the sparsity constrained optimization with the help of our new merit function.}}
\end{itemize}

\subsection{Organization} The rest of the paper is organized as follows. In the next section, we introduce some basic concepts including subdifferential, the generalized Hessian and P-matrix.  \Cref{sec:va} presents the calculations of the gradient and generalized Hessian of the merit function $f_r$ and also establishes the relationship between a solution  to the LCP and a stationary point of $f_r$. We prove several properties of the sparse LCP (\ref{sncp}) via {{the sparsity constrained optimization (\ref{slco})}} in  \Cref{sec:slcp}, including the existence and the boundedness of the solution set to the sparse LCP,  the boundedness of the level set of  $f_r$ over $S$ as well as  the relationship between a solution  to the sparse LCP and {{a stationary point of its sparsity constrained  model}}. In \Cref{sec:nhtp},  we recall the method \NHTP\  and establish its convergence results. Extensive numerical experiments of \NHTP\ solving sparsity constrained models and some concluding remarks are given in the last two sections.

\subsection{Notation} 
{{We  end this section with some notation to be employed throughout the paper.}}   Let $\Diag(x)$ be the diagonal matrix with diagonal elements being from $x$. {{Given two vectors $x,z\in\R^n$, we have the following notation}}
 \begin{eqnarray}\label{notation}
 &~~&\begin{array}{rllrll}
  \mathbb N&:=&\{1,2,\cdots,n\},& \supp(x)&:=&\{i\in\mathbb N:~x_i\neq 0\}\\
 |x|&:=& (|x_1|,\cdots,|x_n|)^\top, &x \circ z&:=& ( x_1z_1,\cdots, x_nz_n)^\top, \\
x_-&:=&\left[(x_1)_-,\cdots,(x_n)_-\right]^\top,&  x^r_+ &:=&\left[((x_1)_+)^r,\cdots,((x_n)_+)^r\right]^\top.
 \end{array}
 \end{eqnarray} 
Note that $x^r_+=(x_+)^r$. For a set $T$, its complementary set is $T^c$ and cardinality is $|T|$.  Denote  $M_{T}$ as the sub-matrix containing the columns of $M$ indexed on ${T}$ and $x_T$ as the sub-vector containing elements of $x$ indexed on ${T}$. However, $M_i$ represents the $i$th row of $M$. In addition, let $e_i$ be the vector with $i$th element being one and {{remaining}} elements being zeros and $e$ be the vector with all elements being ones. Furthermore, write $M_{T_1,T_2}$ as the sub-matrix containing the rows of $M$ indexed on ${T_1}$ and  columns of $M$ indexed on ${T_2}$. Write $M_{T}^\top:=(M_{T})^\top$ and $M_{T_1,T_2}^\top:=(M_{T_1,T_2})^\top$, the transpose of $M_{T}$ and $M_{T_1,T_2}$, respectively.   In particular,  $\nabla_{T} f(x):=(\nabla f(x))_{T}$ and $\nabla^2_{TT} f(x):=(\nabla^2 f(x))_{TT}$, where $\nabla f(x)$ and $\nabla^2 f(x)$ are the gradient and Hessian of $f(x)$. {{Given a matrix $M$, ${\rm rank}(M)$ is the rank and  $M\succeq 0$ (resp. $M\succ 0$) means it is
positive semidefinite (resp. definite). Particularly, we write $A\succeq B$ if $A-B \succeq 0$.  Finally,  define a set $\Xi(\cdot,\cdot)$  by
\begin{eqnarray}\label{set-xi}
\Xi(a,b):=\left\{\begin{array}{rl}
\{b^2_+\},&a>0,\\
{\rm co}\{1,b^2_+\},&a=0,\\
\{1\},&a<0,
\end{array}
\right. 
 \end{eqnarray}
where ${\rm co}\Omega$ is the convex hull of $\Omega$.  Note that  $\Xi(a,b)\neq \Xi(b,a)$ generally.}}

\section{Preliminaries}\label{sec:pre}
In order to analyse functions $\phi_r$ and $f_r$, we first introduce the concept of lower semi-continuity \cite[Definition 4.2]{MN}.
An extended-real-valued function $\varphi:\R^n\rightarrow \R \cup \{+\infty\}$ is lower semi-continuous (l.s.c.) at $\overline{x}\in\R^n $
 if for every $\epsilon\in\R$ with $\varphi(\overline{x})>\epsilon$, there is $\delta>0$ such that
$$\varphi(x)>\epsilon~~~~~{\rm~ for~all}~ x\in U(\overline{x}, \delta):=\{x\in\R^n: \|x-\overline{x}\|<\delta\}.$$
We simply say that $\varphi$ is lower semi-continuous if it is l.s.c. at every point of $\R^n$. From  \cite[Definition 8.3]{RW1998}, for a  proper and l.s.c. function {{$\varphi:\R^n\rightarrow \R$,}} the regular subdifferential and the limiting subdifferential  are  respectively defined as
 \begin{eqnarray*}
   \widehat{\partial}\varphi(x)&= & \left\{v\in {{\R^n:}}~~
 \underset{z(\neq x)\rightarrow x}{\rm liminf} \dfrac{\varphi(z)-\varphi(x)-\langle v, z-x\rangle
 }{\|z-x\|}\geq 0\right\},  \\
   \partial \varphi(x)&= & \underset{z\stackrel{\varphi}{\rightarrow} x}{\limsup} ~\widehat{\partial}\varphi(z)= \left\{v\in  {{\R^n}}: \begin{array}{r}
                              \exists~z\stackrel{\varphi}{\rightarrow} x,~v_j\in \widehat{\partial}\varphi(z_j)~\text{with}~ v_j\rightarrow v
                           \end{array} \right\},
 \end{eqnarray*}
{{where
$z\stackrel{\varphi}{\rightarrow} x$ means both $z\rightarrow x$  and $\varphi(z)\rightarrow \varphi(x)$. If $\varphi$ is convex, then the limiting subdifferential  is also known to be a subgradient.  If it is continuously differentiable, then the limiting subdifferential  is also known as the gradient, i.e., $\partial \varphi(x)=\{\nabla \varphi(x)\}$.}}

%
The next concept is the  (Clarke) generalized Jacobian or the generalized Hessian.  Consider a locally Lipschitz function $F : \R^n\rightarrow \R^m$
  and fix $x\in\R^n$. The  generalized Jacobian \cite{C75} of $F$ at $x$ is the following set of
$m\times n$ matrices:
\begin{eqnarray}
\label{gen-jac} \partial F(x) &=& {\rm co} \left\{\lim \nabla F(x^k) : x^k \rightarrow x, x^k \in D_F \right\}, \end{eqnarray}
where $\nabla F(x^k)$ stands for the classical Jacobian matrix of $F$ at $x^k$ and $D_F$ denotes the set of all the points where $F$ is differentiable.  The generalized Hessian \cite[Definition 2.1]{HSN84} of a continuously differentiable function $\varphi$ at $x$ is defined by
$$\partial^2 \varphi(x):=  \partial (\nabla \varphi(x)).$$
As stated in  \cite[Example 2.2]{HSN84},  $\varphi$ is convex on $\Omega$ if and only if  $\partial^2 \varphi(x)$ is positive semidefinite for all $x\in\Omega$. Here, $\partial^2 \varphi(x)$ is positive semidefinite at $x$ if all elements in $\partial^2 \varphi(x)$ are positive semidefinite. Now we are ready to give our first result with regard  to the first and second order information of functions $a_+^r$ and $|a_-|^r$.
\begin{proposition}\label{pro-phi} {{The following results}} hold for functions $a_+^r$ and $|a_-|^r$.
\begin{itemize}
\item[1)]For any $r>2$, both $a_+^r$ and $|a_-|^r$ are twice continuously differentiable and
\begin{eqnarray*}
\begin{array}{llllll}
\nabla (a_+^r)&=&r a_+^{r-1}, &\nabla^2 (a_+^r)&=&r(r-1) a_+^{r-2},\\
\nabla (|a_-|^r)&=&-r |a_-|^{r-1},& \nabla^2 (|a_-|^r)&=&r(r-1) |a_-|^{r-2}.
\end{array}
\end{eqnarray*}
 \item[2)]For $r=2$, both $a_+^2$ and $|a_-|^2$ are  continuously differentiable and
\begin{eqnarray*}
\begin{array}{llllll}
\nabla (a_+^2)&=&2 a_+ , &\partial^2 (a_+^2)&=&\left\{
\begin{array}{ll}
\{2a_+/a\},&a\neq 0,\\
\left[0,2\right],&a=0,
\end{array}
\right. \\
\nabla (|a_-|^2)&=&2 a_-,&\partial^2 (|a_-|^2)&=& \left\{
\begin{array}{ll}
\{2a_-/a\},&a\neq 0,\\
\left[0,2\right],&a=0.
\end{array}
\right.
\end{array}
\end{eqnarray*}
\end{itemize}
\end{proposition}
Based on above results, we have {{ the following properties}} of $\phi_r$.
\begin{proposition}\label{pro-phi-1} {{The following results}} hold for $\phi_r$ defined by (\ref{func-phi}).
\begin{itemize}
\item[1)]For  any $r\geq2$, $\phi_r$ is continuously differentiable on $ \mathbb{R}\times \mathbb{R}$ with  
    \begin{eqnarray*}\label{grd-phi}
\nabla \phi_r(a,b)=\left[
\begin{array}{cc}
a_+^{r-1}b_+^r-|a_-|^{r-1}\\
a_+^rb_+^{r-1}-|b_-|^{r-1}
\end{array} \right].
\end{eqnarray*}
In addition, $\nabla \phi_r(a,b)=0$ if and only if $\phi_r(a,b)=0$.
\item[2)] For  any $r>2$, $\phi_r$ is twice continuously differentiable on $ \mathbb{R}\times \mathbb{R}$   with 
\begin{eqnarray*}\label{hess-phi}
\nabla^2\phi_r(a,b)=\left[
\begin{array}{cc}
(r-1)(a_+^{r-2}b_+^r+|a_-|^{r-2})&r  a_+^{r-1}b_+^{r-1} \\
r  a_+^{r-1}b_+^{r-1}& (r-1)(a_+^rb_+^{r-2}+|b_-|^{r-2})
\end{array} \right].
\end{eqnarray*}

\item[3)] {{For  $r=2$, the generalized Hessian of $\phi_2(a,b)$  at $(a,b)\in \mathbb{R}\times \mathbb{R}$  has the form  
\begin{eqnarray*}\label{hess-phi2}
\partial\nabla\phi_r(a,b)\subseteq\left\{
\left[
\begin{array}{cc}
u&2a_+b_+\\
2a_+b_+&v
\end{array} \right]: u\in\Xi(a,b), v\in\Xi(b,a)
\right\},
\end{eqnarray*}
where $\Xi(\cdot,\cdot)$ is defined as (\ref{set-xi}).}}
\end{itemize}
\end{proposition}
 {{The proofs of the above two propositions are omitted since they are quite simple. Now, we   compare $\phi_r$ with some other famous NCP  functions.}}
\begin{remark}\label{remark1} We summarize several types of  NCP functions  as follows.
\begin{itemize}
  \item[i)] FB-type functions:
  \begin{eqnarray*}
  \begin{array}{llll}
 \phi_{FB}(a,b)&:=&\left[a^2+b^2\right]^{1/2}-a-b,&\\
   \phi_{FB}^\nu (a,b)&:=&\left[(a-b)^2+\nu ab\right]^{1/2}-a-b,&\nu \in(0,4),\\
   \phi_{FB}^\theta(a,b)&:=&\left[\theta(a-b)^2+(1-\theta) (a+b)^2\right]^{1/2}-a-b,&\theta\in[0,1],\\
   \phi_{FB}^\kappa(a,b)&:=&\left[a^\kappa+b ^\kappa\right]^{1/\kappa}-a-b,&\kappa>1.
  \end{array}
  \end{eqnarray*}
  More details of the above functions can be found in \cite{F92,KK98,LW09} and \cite{CP08}, respectively. The most well-known function among them is the Fischer-Burmeister function $\phi_{FB}$. It was first introduced by Fischer in \cite{F92} and widely used in designing semismooth Newton-type methods for solving mathematical programming with complementarity conditions.   All those functions share a similar {{mathematical formula and hence enjoy similar properties. At the origin, they are nondifferentiable and  have unbounded Hessian.}}
  \item[ii)]Natural residual (minimum function) \cite{CH97}: $$\phi_{\min}(a,b):=2\min\{a,b\}=a+b- [(a-b)^2 ]^{1/2}.$$ This function is simple  but  contains little of the second order information. It is differentiable everywhere {{except at the origin and along the line $a=b$}}.
  \item[iii)]A convex combination function \cite{CCK00}:  $$\phi^{\lambda}(a,b):=\lambda \phi_{FB}(a,b)+(1-\lambda)a_+b_+$$ with $\lambda\in(0,1)$. It is nondifferentiable at $\{(a,b): a\geq0,b\geq0,ab=0\}$.
  \item[iv)] A function proposed in \cite{MS93}:
  \begin{eqnarray*}\label{func-phi-ms93}
\phi^\alpha(a,b)&:=&(ab)^2+\alpha\max\{0,-a,-b \}^{2},
\end{eqnarray*}
where $\alpha>0$. It is continuously differentiable everywhere.
  \item[v)]A class of functions {{proposed in}}  \cite{KYF97},
\begin{eqnarray*}\label{func-phi-kyf} \phi^p_I(a,b)&:=&(ab)_+^p+\left[|a_-|+|b_-| \right]^{p},\\
\phi^p_{II}(a,b)&:=&(ab)_+^p+\left[|a_-|^2+|b_-|^2 \right]^{p/2},\\
\phi^p_{FB}(a,b)&:=&(ab)_+^p+\left[\phi_{FB}(-a,-b)\right]_+^{p},\\
\phi^p_{\max}(a,b)&:=&(ab)_+^p+\max\{0,-a,-b \}^{p},
\end{eqnarray*}
where $p>1$, which is continuously differentiable up to $(p-1)$th order.
\end{itemize}
When these functions in i)-iv) are applied  to deal with the linear/nonlinear complementarity problems, their squared version $\phi^2$ are used and thus are  continuously differentiable everywhere but not twice continuously differentiable.   Compared with those functions, $\phi_r$ defined as (\ref{func-phi}) is  also continuously differentiable for any $r\geq2$ as well as twice continuously differentiable everywhere  for any $r>2$. Moreover, it has bounded Hessian near the origin. Compared with those functions in v), $\phi_r$ has a different first term $a_+^rb_+^r$ and removes the crossed term  $|a_-||b_-|$.  This allows calculations of first and second order derivatives of $\phi_r$ easier. Note that the crossed term can be gotten rid of in $\phi^p_{II}$ only when $p=2$ and in $\phi^p_{I}, \phi^p_{\max}$ only when $p=1$. More interestingly, when the linear mapping $M$  is  positive semi-definite, $\phi_r$ enables $f_r$ to be convex, {{see 4)}} in \Cref{pro-LCP}, which means $\min_x f_r(x)$ is an unconstrained convex optimization with the objective function being  continuously differentiable. 

In addition,  similar to  (\ref{slco}) with merit function $f_r$ being created by $\phi_r$, we can derive different sparsity constrained models with merit functions being constructed by different NCP functions. However, numerical experiments (see \Cref{subsec:ncp}) show that the model with our new merit function $f_r$ outperforms the others.
\end{remark}
 
To end this section, {{we recall the concepts}} of the P-matrix, P$_s$-matrix and Z-matrix,   which play an essential role in subsequent analysis.
\begin{definition}\label{P-matrix} Let $s\leq n$ be a  given integer. A matrix $A\in\R^{n\times n}$ is
\begin{itemize}
  \item[1)]a P-matrix if all of its principal minors are positive \cite{FP62}.
  \item[2)] {{a P$_s$-matrix if all of its principal minor of order up to $s$ are positive.}}
  \item[3)]a Z-matrix if its off-diagonal elements are non-positive \cite{CPS92}.
\end{itemize}
\end{definition}
 {{If $A$ is a  P-matrix, then so are each of its principal sub-matrices and their transpose. Also, a P-matrix  must be a P$_s$-matrix, but not vice versa. 
The equivalent expression of P/P$_s$-matrix is stated below.}}
\begin{proposition}\label{P-matrix}Let $s\leq n$ be a  given integer. A matrix $A\in\R^{n\times n}$ is
\begin{itemize}
  \item[1)]a P-matrix if and only if,  for each nonzero $x\in\R^n$,  there is an index $i$ such that $x_i(Ax)_i>0$.
  \item[2)]a P$_s$-matrix if and only if,  for each nonzero $x\in\R^n$  {{with $\|x\|_0\leq s$,}}  there is an index $i$ such that $x_i(Ax)_i>0$.
\end{itemize}
\end{proposition}
\section{Variational analysis}\label{sec:va} The first issue that we confront is the differentiability of $f_r$, therefore, we start with  calculating its gradient and (generalized) Hessian.

\subsection{Subdifferentials' calculation}
\Cref{pro-phi} and \Cref{pro-phi-1} enable  us to claim the following   proposition regarding  the first and second order information of $f_r$ in (\ref{sco-obj}). Hereafter, for notational simplicity, we denote $ y:=Mx+q$.
\begin{lemma}\label{pro-LCP} For $f_r$ as in (\ref{sco-obj}),  {{the following results}} hold.
\begin{itemize}
\item[1)] For any $r\geq2$, $f_r(x)$ is continuously differentiable  with  
\begin{eqnarray}\label{grd-phi-lcp}
&&\nabla f_r(x)=  x_+^{r-1}\circ y_+^r - |x_-|^{r-1} + M^\top \left[x_+^{r}\circ y_+^{r-1} - |y_-|^{r-1}\right].
\end{eqnarray}
\item[2)] For any  $r>2$, $f_r(x)$ is twice continuously differentiable with  
\begin{eqnarray}\label{hess-phi-lcp}
\nabla^2f_r (x)&=&r\left[ \Diag(x_+^{r-1}\circ y_+^{r-1}) M+M^\top\Diag(x_+^{r-1}\circ y_+^{r-1} )\right]\\
& +&(r-1) \Diag\left(x_+^{r-2}\circ y_+^{r}+|x_- |^{r-2}\right)   \nonumber\\
& +&(r-1)M^\top \Diag\left(x_+^{r}\circ y_+^{r-2}+|y_- |^{r-2}\right)  M.\nonumber
\end{eqnarray}
\item[3)] For  $r=2$, the generalized Hessian $\partial^2 f_2(x)$ takes the form
\begin{eqnarray}
\label{hess-phi2-lcp}\partial^2f_2 (x)&\subseteq&\Big\{ 2\left[ \Diag(x_+\circ y_+) M+M^\top\Diag(x_+\circ y_+ )\right] +\\
&& ~~\Diag (\xi ) + M^\top \Diag (\zeta )  M:~~\xi\in\Omega_\xi(x),~ \zeta\in\Omega_\zeta(x)\Big\},\nonumber
\end{eqnarray}
where $\Omega_\xi(x)$  and $\Omega_\zeta(x)$ are given by {{
\begin{eqnarray}
\label{xilcp}
\Omega_\xi(x)&:=&\left\{\xi\in\R^n:~\xi_i\in \Xi(x_i,y_i)\right\},\\
\label{zetalcp}\Omega_\zeta(x)&:=&\left\{\zeta\in\R^n:~\zeta_i\in \Xi(y_i,x_i)\right\},
\end{eqnarray}
where $\Xi(\cdot,\cdot)$ is defined as (\ref{set-xi}).}}
 \item[4)]For any $r\geq2$, $f_r(x)$ is convex if $M$ is positive semidefinite.

\end{itemize}
\end{lemma}

\subsection{Stationary points}
This subsection reveals relationship between the solutions to the LCP and the stationary points of $f_r$. We say a point $x^*$ is a stationary point of $f_r$ if it satisfies
\begin{eqnarray}\label{s-set}x^*\in \{x\in\R^n:~\nabla f_r(x)=0 \}=:\mathcal{G}_f.\end{eqnarray}
Moreover, we  say the LCP  is feasible if
\begin{eqnarray}\label{CP-fea}
\fea:= \{x\in\R^n:~x\geq0,\  Mx+q \geq0  \}\neq \emptyset.
\end{eqnarray}
Based on \cite[Proposition 3.1.5]{CPS92}, the LCP  is feasible  for all $q\in\R^n$ if and only if there is an $x$ such that $x>0, Mx>0$. According to \cite[Definition 3.1.4]{CPS92}, the matrix satisfying such condition is called S-matrix. One could easily derive that if $ \sol \neq \emptyset$, then
\begin{eqnarray}\label{sol-sco} \sol ={\rm argmin}_{x}  \ f_r(x).\end{eqnarray}
Because of this, it is obvious that $\sol\subseteq \mathcal{G}_f$ since an optimal solution is also a stationary point,  {{while the converse is not true in general. However, under some assumptions, we can claim that these two sets coincide.}}
\begin{theorem} \label{the-sta}For any given $q\in\R^n$, we have the following results.
\begin{itemize}
\item[1)]  If $M$ is positive semidefinite and  \emph{\texttt{fea}}$(M,q)$ is nonempty, then \emph{\texttt{sol}}$(M,q)=\mathcal{G}_f$ is nonempty as well.

\item[2)] If $M$ is a P-matrix, then \emph{\texttt{sol}}$(M,q)=\mathcal{G}_f=\{x^*\}$, where $x^*$ is the unique solution to \emph{\texttt{sol}}$(M,q)$. 

\end{itemize}
 \end{theorem}
\section{Sparse LCP}\label{sec:slcp}
Now we center on the sparse LCP (\ref{sncp}) and its corresponding sparsity constrained optimization  (\ref{slco}) through the proposed merit function $f_r$. We start studying  the existence and boundedness of the solution set to the sparse LCP. Hereafter, we say the sparse LCP is feasible if
\begin{eqnarray}\label{sCP-fea}
&&~~\texttt{fea}_s(M,q):=\texttt{fea}(M,q)\cap S=\left\{x\in\R^n:~x\geq0,\  Mx+q \geq0, \ \|x\|_0\leq s \right\}
\end{eqnarray}
is nonempty. One can see that, for example, if $M$ is a  matrix with  {{all entries being positive}}, then the sparse LCP is feasible for any $q\in\R^n$. In fact, for any $x\geq0$ with $\|x\|_0\leq s$, one can find a proper large $\delta$ such that $M (\delta x)+q\geq 0$, which means $\delta x\in\texttt{fea}_s(M,q)$. Some other types of matrices may also guarantee  the feasibility of the sparse LCP. However, we will not explore them in this paper and simply assume that \texttt{fea}$_s(M,q)$ is nonempty   {{in the sequel}}. 

\begin{lemma}\label{exist-s-quad} If \emph{\texttt{fea}}$_s(M,q)$ is nonempty, then  so is
\begin{eqnarray}\label{s-quad}
Q_s(M,q):={\rm argmin}_x~ \langle x, Mx+q \rangle,~~{\rm s.t.} ~ x\in \emph{\texttt{fea}}_s(M,q).
\end{eqnarray}
\end{lemma}
 
\subsection{Existence and boundedness}
 Our first result is about the existence of solutions to the sparse LCP under some assumptions. Note that if $q\geq0$, then $0\in\sol$, a trivial solution. 	In other words, if there is  an $i$ such that $q_i<0$,  then $0\notin\sol$. 
For a point $x$, denote two sets
 \begin{eqnarray}\label{TG}
T:=\supp(x),~~~~\Gamma:=\{i\in\N: M_ix+q_i=0\}.
  \end{eqnarray}
Here, $T$ and $\Gamma$ are depended on $x$. We drop their dependence for notational simplicity. Now, we give the results about the existence of a solution to the sparse LCP.
 \begin{theorem}
  \label{existence-main}  Assume \emph{\texttt{fea}}$_s(M,q)$ is nonempty, which means there exists an $x\in Q_s(M,q)$.  Then $x\in \emph{\texttt{sol}}(M,q)\cap S$ if one of   {{the following conditions}} holds
\begin{itemize}
 \item[1)] $\|x\|_0=s$, $M$ is a symmetric Z-matrix with ${\rm rank} (M_T)=|T|$ and $q_T\leq0$.
  \item[2)] $\|x\|_0<s$, $M$ is a symmetric Z-matrix with ${\rm rank}(M_{T \Gamma})=|\Gamma|$ and $q_T\leq0$.
  \item[3)] $\|x\|_0<s$, $M$ is positive semidefinite with  ${\rm rank}(M_{T \Gamma})=|\Gamma|$.
\end{itemize}
\end{theorem}
{{It is worth mentioning that in \Cref{existence-main} 2), $T\subseteq\Gamma$ from the proof of 2) in \Cref{proof-exist}, while the assumption that $M_{T\Gamma }$ has full column rank requires  $|T| \geq |\Gamma|$. Therefore, there is $T =\Gamma$.}} Next result exhibits another sufficient condition to guarantee the existence of a solution to the sparse LCP.
\begin{theorem}
  \label{existence}
 Assume $M$ is a $P_s$-matrix with all entries being nonnegative.  If $|\theta|\leq s$, where $\theta:=\{i\in\N:q_i<0\}$,  then $\emph{\texttt{sol}}(M,q) \cap S$ is nonempty and contains a unique $x^*$ such that $\supp(x^*)\subseteq \theta$.
\end{theorem}
We now have the boundedness of  the following level set. This suffices to show the boundedness of the solution set $(\sol \cap S)$ to (\ref{sncp}).
  \begin{theorem}\label{pro-bound-s} If $M$ is a  P$_s$ matrix, then the  level set 
 \begin{eqnarray} \label{Ls}\mathcal{L}_s(f_r,\gamma):=\{x\in S:~f_r(x)\leq\gamma\}  \end{eqnarray}
 is bounded for any $\gamma\geq0$. Moreover, $(\emph{\texttt{sol}}(M,q) \cap S)\subseteq \underset{x\in S}{\rm argmin} f_r(x)$ are both bounded.
 \end{theorem}


\subsection{Optimality Conditions}
  \Cref{pro-bound-s} indicates an optimal solution of   (\ref{slco})  must exist if $M$ is a  P$_s$ matrix. In addition,
it follows from \cite[Theorem 2.8]{PXZ15} that an optimal solution $x^*\in S$  of   (\ref{slco}) satisfies  \begin{eqnarray}\label{necc-opt-cod}-\nabla f_r(x^*)\in N_{S}(x^*), \end{eqnarray}
where $N_{S}(x^*)$ is the Bouligand normal cone of $S$ at $x^*$. {{Hereafter, let 
\begin{eqnarray}\label{T*supp} T_*:=\supp(x^*)\end{eqnarray} for notational convenience. From \cite[Table 1]{PXZ15},  the condition \eqref{necc-opt-cod}}} is equivalent to
\begin{eqnarray}\label{necc-opt-cod-1}
&&\nabla_{i} f_r(x^*)\left\{\begin{array}{cc}
  =0, & i\in T_*, \\
  \in\R, & i\notin T_*,
\end{array}\right. {\rm if}~ \|x^*\|_0=s\  \ {\rm and}\  \ \nabla f_r(x^*)=0\ {\rm if}~ \|x^*\|_0<s.
  \end{eqnarray}
We call a point a stationary point of (\ref{slco}) if it satisfies (\ref{necc-opt-cod-1}). The next theorem reveals the relationship  between a stationary point and a solution to (\ref{sncp}).
  \begin{theorem} \label{the-sta-s} A solution to (\ref{sncp}) is also a stationary point of (\ref{slco}). Conversely, assume that $M$ is a Z-matrix.  Then a stationary  point $x$  of (\ref{slco})  is also a solution to (\ref{sncp}) if  there is a nonzero vector $v\in\R^{|T_+|}$ such that $M_{\Gamma_+^cT_+}v\geq0$ and $M_{T_+T_+}$ is  positive semidefinite, where $T_+:=\{i\in\N:~x_i>0\}$ and $\Gamma_+:=\{i\in\N:~M_ix+q_i>0\}$.
 \end{theorem}
\begin{remark}
With regard to {{the above theorem}}, some comments can be made.
\begin{itemize}
\item[i)] If $\Gamma_+^c\subseteq T_+$ in  \Cref{the-sta-s}, then $M_{T_+T_+}$ being  positive semidefinite indicates    that there always exists a nonzero vector $v\in\R^{|T_+|}$ such that $M_{\Gamma_+^cT_+}v\geq0$.
\item[ii)] {{We give some explanations about $T,T_+,\Gamma_+$ and $\Gamma$. Let \begin{eqnarray*}\label{sets--}
T_-&:=&\{i\in\N:~x_i<0\}, \hspace{5mm}\Gamma_-~:=\{i\in\N:~M_ix+q_i<0\}.
\end{eqnarray*}
Then $T_+$ and $T_-$ capture the indices of positive and negative elements of $\bx$, and hence $T_+\cup T_-=T$ by \eqref{TG}. While $\Gamma+$ and $\Gamma_-$ contain the indices of positive and negative elements of $M\bx+q$, and hence $\Gamma_+\cup \Gamma_-=\Gamma^c$ by \eqref{TG}.}}
\item[iii)] If a stationary point $x$  of (\ref{slco}) satisfies $Mx+q\geq0,$ then $\Gamma_{-}=\emptyset$. This together with (\ref{g-a-b-t}), (i.e., $|x_{T_-}|^{r-1}=0$ leading to $\bx\geq 0$) suffices to show that $x$ is also a solution to (\ref{sncp}) if  $M_{T_+T_+}$ is  positive semidefinite. As a consequence, the other assumptions can be neglected.
\end{itemize}
\end{remark}
 We end this section with establishing the relationship between a stationary point and a local/global solution to (\ref{slco}) by the following theorem.
 \begin{theorem}\label{theorem-local-stationary}
 Assume that $M$ is positive semidefinite. Consider a point $x^*\in S$.
 \begin{itemize}
   \item[1)]If $\|x^*\|_0<s$, then it is a stationary point if and only if it is a globally optimal solution to (\ref{slco}). If  {{we further assume that \emph{\texttt{fea}}$(M,q)$ is nonempty, then the stationary point  satisfies $x^*\in(\emph{\texttt{sol}}(M,q)\cap S).$}}
   \item[2)]If $\|x^*\|_0=s$, then it is a stationary point if and only if it is a locally optimal solution to (\ref{slco}). If  {{we further assume that $M_{T_*T_*}$ is nonsingular, then the stationary point $x^*$ is a unique optimal solution to (\ref{slco}) with $r=2$ on $\R_{T_*}:=\{x\in\R^n:\supp(x)\subseteq T_*\}.$}}
 \end{itemize}
  \end{theorem}

 \section{Newton Hard-Thresholding Pursuit}\label{sec:nhtp}
 We now turn our attention to the solution method, Newton Hard-Thresholding Pursuit (\NHTP), for (\ref{slco}). The method is adopted from \cite{ZXQ19}. To implement the method, we first define some notation. \begin{eqnarray}\label{Tu}
 ~~~~~~\T(x,\eta):= \left\{ T\subseteq \N:\begin{array}{l}
T~ \text{contains the indices of}~s~\text{largest elements of}~|z|\\
|T|=s, ~ \text{where}~ z:=x-\eta \nabla f_r(x) 
\end{array} \right\},\end{eqnarray}
where $\eta>0$. Note that $T$ may not be unique since the $s$th largest element of $|z|$ might be multiple.  For any given $T \in \T(x; \eta)$, we define a nonlinear equation:
\begin{eqnarray}\label{station}
F_{\eta}(x;T):=
\left[\begin{array}{c}
\nabla_{ T}  f_r (x)\\
x_{ T^c }\\
\end{array}\right]=
0.
\end{eqnarray}
One advantage of defining the function $F_\eta(x; T)$ is that if a point $x$ satisfies $F_\eta(x; T)=0$ for a given $T$ then  it  satisfies (\ref{necc-opt-cod-1}), a stationary point. In addition, this is an  equation system that allows us to perform the Newton method. 
\subsection{Framework of \NHTP}
Suppose $x^k$ is the current approximation to a solution of (\ref{station}) and $T_k$ is chosen
from $\T(x^k; \eta)$.
Then Newton's method for the  equation (\ref{station}) takes the following form to get the direction $d^{k}$:
\begin{equation} \label{Newton-Method}
\nabla F_{\eta} (x^k; T_k) d^k  = - F_{\eta} (x^k; T_k),
\end{equation}
where $\nabla F_{\eta} (x^k; T_k)$ is the Jacobian of $F_{\eta}(x; T_k)$ at $x^k$ and admits the following form:
\begin{equation} \label{Jacobian-Matrix}
\nabla F_{\eta} (x^k; T_k) = \left[
 \begin{array}{cc}
  \nabla^2_{T_kT_k} f_r(x^k) \ & \ \nabla^2_{T_k T_k^c} f_r(x^k)  \\ [0.6ex]
  0 & I_{n-s}
 \end{array}
\right],
\end{equation}
and $\nabla^2  f_r(x)$ is the Hessian of $f_r(x)$ when $r>2$ and   a matrix  from the generalized Hessian $\partial^2 f_2(x)$ when $r=2$. It is worth mentioning that the choice of  $\nabla^2  f_2(x^k)$  does not affect the method proposed in \Cref{alg:nhtp} and its convergence results.  Substituting (\ref{Jacobian-Matrix}) into
(\ref{Newton-Method}) yields
\begin{eqnarray}\label{sequation-k-0}\left\{
\begin{array}{rcl}
\nabla_{{T_{k}} {T_{k}}}^{2} f_r (x^k) d^k_{T_{k}} &=& \nabla_{   T_{k}  T^c_{k} }^{2} f_r (x^k)x^{k}_{  T^c_{k}}-\nabla_{ T_{k}}  f_r (x^k), \\
d^{k}_{  { T }^c_{k}}&=&-x^{k}_{ { T }^c_{k}}.
\end{array} \right.
\end{eqnarray}
After we get the direction, in order to guarantee the next point $x^{k+1}$ to be feasible, namely, $x^{k+1}\in S$, we update it by  {{using the following scheme:}}
\begin{equation} \label{form4}
x^{k}( \alpha):=\left[\begin{array}{cc}
x^{k}_{T_{k}}+ \alpha d^{k}_{T_{k}}\\
0
\end{array}\right] 
\end{equation}
for some $\alpha\in(0,1]$. Now we summarize the whole framework of \NHTP\ in \Cref{alg:nhtp}.  
 
\begin{algorithm}[H]
\caption{\NHTP:  Newton  Hard-Thresholding Pursuit}
\label{alg:nhtp}
\begin{algorithmic}
\STATE{Initialize $x^0$. Choose $\eta,  \gamma>0, \sigma\in(0,0.5), \beta\in(0,1)$ and $K$. Set $k\Leftarrow0$.}
\WHILE{The halting condition does not hold and $k\leq K$}
\STATE{ \emph{Hard-Thresholding Pursuit:} Choose $T_{k}\in\T(x^k,\eta)$ {{in (\ref{Tu}).}}}
\STATE{ \emph{Descent Direction Search:} Update $d^k$ by solving (\ref{sequation-k-0}) if  it is solvable and}
\STATE{ \parbox{0.7\textwidth}{\begin{equation} \label{form1} \langle \nabla_{T_{k} } f_r(x^{k}),  d^{k}_{T_{k} }\rangle\leq-\gamma \|d^{k}\|^2+ \|x^{k}_{ T^c_{k} }\|^2/(4\eta).\end{equation}} }  
\STATE{  Otherwise, update $d^k$ by}
\STATE{ \parbox{0.7\textwidth}{\begin{equation} \label{form2} d^k_{T_{k}}  = -\nabla_{ T_{k}}  f_r (x^k),~~~~d^{k}_{T^c_{k}}=-x^{k}_{T^c_{k}}.\end{equation}} }
\STATE{ \emph{Step Size Search:} Find the smallest integer $t = 0,1,\ldots$ such that}
\STATE{\parbox{0.7\textwidth}{\begin{equation} \label{form3}f_r(x^{k}(\beta^{t}))\leq f_r(x^k)+ {{\sigma \beta^{t}}} \langle \nabla  f_r(x^{k}),  d^{k} \rangle.\end{equation}} } 
\STATE{ Set $\alpha_k=\beta^{t}$ and update  $x^{k+1}=x^{k}( \alpha_k)$ by (\ref{form4}).}
\ENDWHILE
\RETURN the solution $x^k$.
\end{algorithmic}
\end{algorithm}

Some comments can be made based on \Cref{alg:nhtp}. Note that, because of (\ref{form4}), namely, $x^{k+1}=x^{k}( \alpha_k)$, we always have 
\begin{equation} \label{form5}\supp(x^{k+1})\subseteq T_k.\end{equation} 

{\bf (a) Computational complexity.} In \emph{Hard-Thresholding Pursuit} step,  we only pick $s$ {{indices of $s$ largest elements of $|x^k-\eta \nabla f_r (x^k)|$  }}to form $T_{k}$, which allows us to use  \texttt{mink} function in MATLAB (2017b or later version) whose computational complexity is $\mathcal{O}(n+s\log s)$. In \emph{Descent Direction Search} step, from $\supp(x^{k})\subseteq T_{k-1}$ by (\ref{form5}), the first equation of (\ref{sequation-k-0}) can be rewritten as
\begin{eqnarray}\label{form6}\nabla_{{T_{k}} {T_{k}}}^{2} f_r (x^k) d^k_{T_{k}} 
&=& \nabla_{   T_{k}  J_{k} }^{2} f_r (x^k)x^{k}_{J_{k}}-\nabla_{ T_{k}}  f_r (x^k),
\end{eqnarray}
where $J_{k}:=T_{k-1}\cap T^c_{k}$ and thus $|J_{k}|\leq |T_{k-1}|=s$. So we 
 {{need to calculate}} $\nabla_{ T_{k}}  f_r (x^k)$, $\nabla_{ T_{k}  J_{k} }^{2} f_r (x^k)x^{k}_{J_{k}}$ and a sub-Hessian $\nabla_{T_k,T_k}^{2}f_r(x^k)$. It follows from (\ref{grd-phi-lcp}), (\ref{hess-phi-lcp}) or (\ref{hess-phi2-lcp}) that the most computational expensive calculations in these three terms  are
 $$M^\top_{T_k}|(M_{T_{k-1}}x_{T_{k-1}}^k+q)_-|^{r-1},~ M^\top_{T_k}  \Diag(z^k) (M_{J_k}x^{k}_{J_{k}}),~ M^\top_{T_k}  \Diag(z^k) M_{T_k},$$
where $z:=(x^k)_+^{r}\circ (y^k)_+^{r-2}+|(y^k)_- |^{r-2}$  or  $z\in\Omega_\zeta(x^k)$. Their computational complexities are $\mathcal{O}(ns),\mathcal{O}(ns)$ and $\mathcal{O}(ns^2)$, respectively.  Moreover, to update $d^{k}_{ T_{k}}$, we also {{need to solve}} the linear equation (\ref{form6}) with $s$ equations and $s$ variables, which has computational complexity about $\mathcal{O}(s^\kappa)$, where $\kappa\in(2,3)$. {{Let $\bar{t}$ be}} the smallest integer satisfying (\ref{form3}) and it often takes the value 1. Overall, the whole computational complexity of each step in \Cref{alg:nhtp} is  $\mathcal{O}(ns^2+s^\kappa+\bar{t}ns)$.

{\bf (b) Halting condition.} {{A halting condition used in  \cite{ZXQ19} is to calculate
\begin{equation} \label{Tolerance-Function}
\mbox{Tol}_{\eta}(\bx^k;\; T_k) := \| F_\eta (\bx^k; T_k) \| + \max_{i \in T_k^c} \left( | \nabla_i f_r(\bx^k) | - x_{(s)}^k/\eta, \ 0  \right)_+,
\end{equation}
where $x_{(s)}^k$ is the $s$th largest element of $|\bx^k|$. If a point $\bx^k$ satisfies that $\mbox{Tol}_{\eta}(\bx^k;\; T_k) =0$, then both terms on the right-hand side of \eqref{Tolerance-Function} are zeros, which imply that $\nabla_{T_k} f_r(\bx^k)=0, \bx^k_{T_k^c}=0$ and $\| \nabla_ {T_k^c} f_r(\bx^k) \|_\infty \leq x_{(s)}^k/\eta$. Hence $\supp(x^{k})\subseteq T_{k}$. These derive the first condition in \eqref{necc-opt-cod-1} if $\|x^{k}\|_0=s$ and $\nabla f_r(\bx^k)=0$ in \eqref{necc-opt-cod-1} if $\|x^{k}\|_0<s$ since $x_{(s)}^k=0$ under such case. Namely,  $x^{k}$ is a stationary point of \eqref{slco}. Therefore, we will terminate \NHTP\ if $\mbox{Tol}_{\eta}(\bx^k;\; T_k)<\texttt{tol}$ in our numerical experiments, where $\texttt{tol}$ is a tolerance (e.g. $10^{-6}$).}}

\subsection{Convergence analysis}
 {{ As shown in \cite[Theorem 8]{ZXQ19},}} to establish the convergence results, the assumptions are relating to the  boundedness of Hessian and existence of the inverse of the Hessian at the limiting point. We first define a parameter to bound the Hessian under mild condition
\begin{eqnarray}\label{bound-hessian}
C:=\sup_{x\in\mathcal{L}_s(f_r,f_r(0))} \sigma_{\max}(\nabla^2 f_r(x)),\end{eqnarray} 
where $\mathcal{L}_s(f_r,f_r(0))$ is the level set given as (\ref{Ls}) and $\sigma_{\max}(A)$ is the maximum singular value of $A$. The following result {{ shows that}} such $C$ is bounded if $M$ is a  P$_s$ matrix.
\begin{lemma}\label{bound-c} If $M$ is a  P$_s$ matrix, then  $C<+\infty$.
 \end{lemma}
 Denote a parametric point $\mu:=(\eta,  \gamma, \sigma, \beta)$ where $\eta>0,  \gamma>0, \sigma\in(0,0.5), \beta\in(0,1)$. {{Based on the above lemma, we have the following convergence results. }}
\begin{theorem}\label{theorem-convergence} Suppose $M$ is a  P$_s$ matrix and also positive semidefinite. Choose $x^0\in\mathcal{L}_s(f_r,f_r(0))$ with $f_r(x^{0})\leq f_r(0)$. Then there exist some $\mu$ such that the following results hold.
\begin{itemize}
\item[1)] $\{f_r(x^k)\}$ is non-increasing and $\{x^k\}$ is bounded.
\item[2)] Any accumulating point, say $x^*$, of the sequence $\{x^k\}$ is a stationary point of (\ref{slco}) and thus a local minimizer by \Cref{theorem-local-stationary}.
\item[3)]  {{If further assume that $\nabla^2_{T_\infty T_\infty} f_r(x^*)$ is invertible for any $T_\infty\supseteq\supp(x^*)$ and $|T_\infty|=s$, then the whole sequence converges to $x^*$ and the Newton direction is always admitted for sufficiently large $k$.}}
\end{itemize}
\end{theorem} 
 
\begin{remark} {{We give some explanations about the conditions in \Cref{theorem-convergence}.}}
\begin{itemize}
\item[i)] {{If $x^*$ is a solution to the sparse LCP, then $\nabla^2 f_r(x^*)=0$ for any $r>2$ by (\ref{hess-phi-lcp}) and $\nabla^2_{T_\infty T_\infty} f_2(x^*)\succeq M_{T_\infty}^\top \Diag(\varsigma) M_{T_\infty}$ for $r=2$ by (\ref{hess-phi2-lcp}). Therefore, the assumption that $\nabla^2_{T_\infty T_\infty} f_r(x^*)$ being invertible for any $T_\infty\supseteq\supp(x^*)$ and $|T_\infty|=s$ does not hold for $r>2$ but holds for $r=2$ most likely.}} This might be a reason that the sparsity constrained model with $f_2$ outperforms the other models with $f_r$ for $r>2$, see \Cref{subsec:r}.
\item[ii)] {{The choice of}} $x^0\in\mathcal{L}_s(f_r,f_r(0))$ with $f_r(x^{0})\leq f_r(0)$ in \Cref{theorem-convergence} is easy to be satisfied. One could choose $x^0=0$ for simplicity. This choice also gives us an initial point when we implement  \Cref{alg:nhtp} {{in the next section.}}
\item[iii)] {{The choices of $\mu$ can be found in \cite{ZXQ19}. More precisely, $\sigma\in(0,1/2), \beta\in(0,1),$ 
$$0<\gamma\leq\min\{1,2C\},~~~~ 0<\eta\leq\min\left\{ \gamma c\beta/C^2, c\beta, 1/(4C)\right\},$$
where $C$ is given by \eqref{bound-hessian} and $c:=\min\{1, \gamma(1-2\sigma)/(C-\sigma\gamma)\}$. Note that those parameters are   dependent on the objective function $f_r$ and $x^0$ (independent of the iterates $x^k, k\geq 1$ and its limit $x^*$). Moreover, the conditions of those parameters are sufficient but not necessary to guarantee the convergence property. Therefore, there is no need to set them to strictly  meet those conditions in practice, not to mention $c$ or $C$ being difficult to calculate. When it comes to the numerical computation, some of them are suggested to be updated iteratively, such as $\gamma_k=10^{-10}$ if $\bx^k_{  T^c_{k}}=0$ and  $10^{-4}$ otherwise.
}}
\end{itemize}
\end{remark}

\section{Numerical Experiments}\label{sec:ne}

In this part, we implement \NHTP\footnote{available at \url{https://github.com/ShenglongZhou/NHTPver2}} {{described in  \Cref{alg:nhtp}   to solve the sparsity}} constrained complementarity problem  (\ref{sncp}). All experiments were conducted by using MATLAB (R2018a) on a desktop of 8GB memory and Inter(R) Core(TM) i5-4570 3.2Ghz CPU.    {{We terminate the proposed method at the $k$th step if it meets one of the following conditions: 1) $\mbox{Tol}_{\eta}(\bx^k;\; T_k) \leq 10^{-6}$, where $\mbox{Tol}_{\eta}(\bx^k;\; T_k)$ is defined as (\ref{Tolerance-Function}); 2) $|f_r(\bx^{k+1})-f_r(\bx^{k})|<10^{-6}(1+|f_r(\bx^{k})|)$ and 3) $k$ reaches the maximum number (e.g., 2000) of iterations.}}  For parameters in \NHTP, we keep all default ones except for $\texttt{pars.eta}$, which is set as $\texttt{pars.eta}=5$ if $n\leq 1000$ and $\texttt{pars.eta}=1$ otherwise for all   numerical experiments.

{{The rest of this section  is organized as follows. We first give four examples to be tested throughout the whole simulations. Since $f_r$ and $S$ in the sparsity constrained model (\ref{slco})  involve parameters $r$ and $s$,   we then run \NHTP\ to see the performance under different choices of $r$ and $s$. Next, we provide two strategies to select a proper $s$ in model \eqref{slco}  in case the sparsity level $s$ is unknown. Followed are the numerical comparisons of   \NHTP\ and two other solvers: half thresholding projection (\HTP) \cite{SZX15} and extra-gradient thresholding algorithm (\ETA) \cite{SZPX15}.  In conclusion, \NHTP\ is capable of producing high quality solutions with fast computational speed when benchmarked against other methods. Finally, to testify the advantage of our new merit function $f_r$, we also apply \NHTP\  to deal with the sparsity constrained model (\ref{slco}) with other merit functions constructed by three existing famous NCP functions: $\phi_{FB}$, $\phi_{\min}$ and $\phi^2_{II}$, see \Cref{remark1}. Numerical comparisons demonstrated that the sparsity constrained model with the new merit function  enables \NHTP\ to run the fastest due to the lowest computational complexity and produce the most accurate solutions.}}

\subsection{Test examples}
Four sparse LCP examples are taken into consideration.  The first three examples have the given `ground truth' sparse solutions $x^*$, while for the last one, the `ground truth' sparse solutions $x^*$ are unknown. It is worth mentioning there are many nonlinear complementarity problems from \cite{MD87,HP90,TFI93, HL02, YHS09, SZX15}, which could be converted to the sparsity constrained optimization through $\phi_r$. We {{had also applied \NHTP\ to}} solve those problems and got the excellent numerical performance. However, we omit the related results to shorten the paper here.

\begin{example}[Z-matrix]\label{z-matrix} Let $M$ and $q$ in  (\ref{sncp}) be given by
 \begin{eqnarray}
  M=I_n-ee^\top/n
  \quad \mbox{and}\quad  q=e/n-e_1,
\nonumber
 \end{eqnarray}
where $I_n$ is the identity matrix of order $n$. Such $M$ is a so-called positive semidefinite Z-matrix and widely used in statistics, which allows that (\ref{sncp}) admits a unique sparse solution $x^*=e_1$  \cite{SZX14}.
\end{example}

\begin{example}[SDP Matrices]\label{sdp-matrix} In  (\ref{sncp}),   a positive semidefinite matrix $M$ and $q$ are given as follows. Let $M=ZZ^\top$ with $Z\in \R^{n\times m}$ whose elements are  generated from the standard normal distribution, where $m\leq n$ (e.g. $m=n/2$).   Then, the `ground truth' sparse solution $x^*$ is produced by the following pseudo Matlab codes:
$$x^* = \texttt{zeros}(n, 1),~\Gamma = \texttt{randperm}(n),~ x^*(\Gamma(1:s)) = 0.1+|\texttt{randn}(s^*,1)|,$$
where $s^*$ is the sparsity level of the solution. We add $0.1$ to generate $x^*$, {{avoiding elements with a tiny scale.}}  Finally, $q$ is obtained by
$$q_i=\left\{\begin{array}{rr}
    -(M x^*)_i, & x_i^*>0, \\
    |(M x^*)_i|, & x_i^*=0.
  \end{array}\right.
$$
 \end{example}

\begin{example}[Nonnegative SDP Matrices]\label{non-sdp-matrix} As stated in  \Cref{existence}, we consider $M$ and $q$ in  (\ref{sncp})  as follows. Let $M=ZZ^\top$ with $Z\in \R^{n\times m}$ whose elements are  generated from the uniform distribution between $[0,1]$, where $m\leq n$ (e.g. $m=n/2$).   Then,  $x^*$ is {{produced as  in}}  \Cref{sdp-matrix} 
and $q$ is obtained by
$$q_i=\left\{\begin{array}{rr}
    -(M x^*)_i, & x_i^*>0, \\
    \texttt{rand}(1), & x_i^*=0.
  \end{array}\right.
$$
 \end{example}

\begin{example}[Nonnegative SDP Matrices without $x^*$]\label{non-sdp-matrix-nox} This example is similar to   \Cref{non-sdp-matrix} but  without given the `ground truth' solution.  Here $M$  is {{generated as  in}}  \Cref{non-sdp-matrix} but with $m=n/4$. Let $\Gamma = \texttt{randperm}(n)$ and $T=\Gamma(1:s^*)$. Then,
 $q$ is obtained by
$$q_i=\left\{\begin{array}{rr}
    -\texttt{rand}(1), & i\in T, \\
    \texttt{rand}(1), & i\notin T.
  \end{array}\right.
$$
 \end{example}

 \subsection{Effect of $r$ with fixing $s=s^*$}\label{subsec:r}
{{The objective function}} $f_r$ involves a parameter $r$. To see the effect of $r$ on (\ref{sncp}), we first compare \NHTP\ solving (\ref{sncp}) under different choices of $r$ but with fixing $s=s^*$ in $S$. Thus, for a given $r$, we write \NHTP\ as \NHTP$_r$. Let $x$ be the solution produced by a method. We say a recovery of this method is successful if $$\|x -x^*\| < 0.01\|x^*\|.$$

For each example, each instance has two deciding factors: $(n, s^*)$. We begin with solving   \Cref{sdp-matrix}  and \Cref{non-sdp-matrix} with fixed $n = 200$  but with increasing sparsity level $s^*$ from $2$ to $44$. For each $(n, s^*)$, we run $500$ independent trials and record the corresponding success rates which is defined by the percentage of the number of successful recoveries {{over all trials.}}

\begin{figure}[H]
\centering
\begin{subfigure}{0.49\textwidth}\centering
  \includegraphics[width=.9\linewidth]{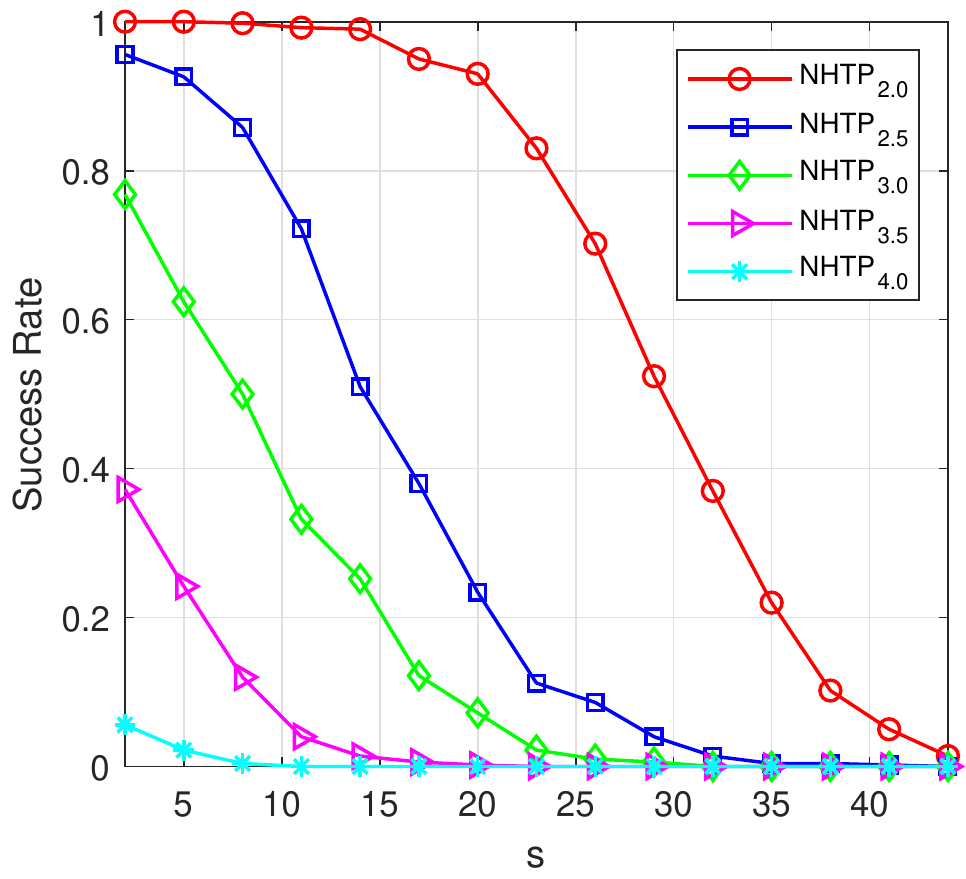}
   \caption{  \Cref{sdp-matrix}}
\end{subfigure}
\begin{subfigure}{0.49\textwidth}\centering
  \includegraphics[width=.9\linewidth]{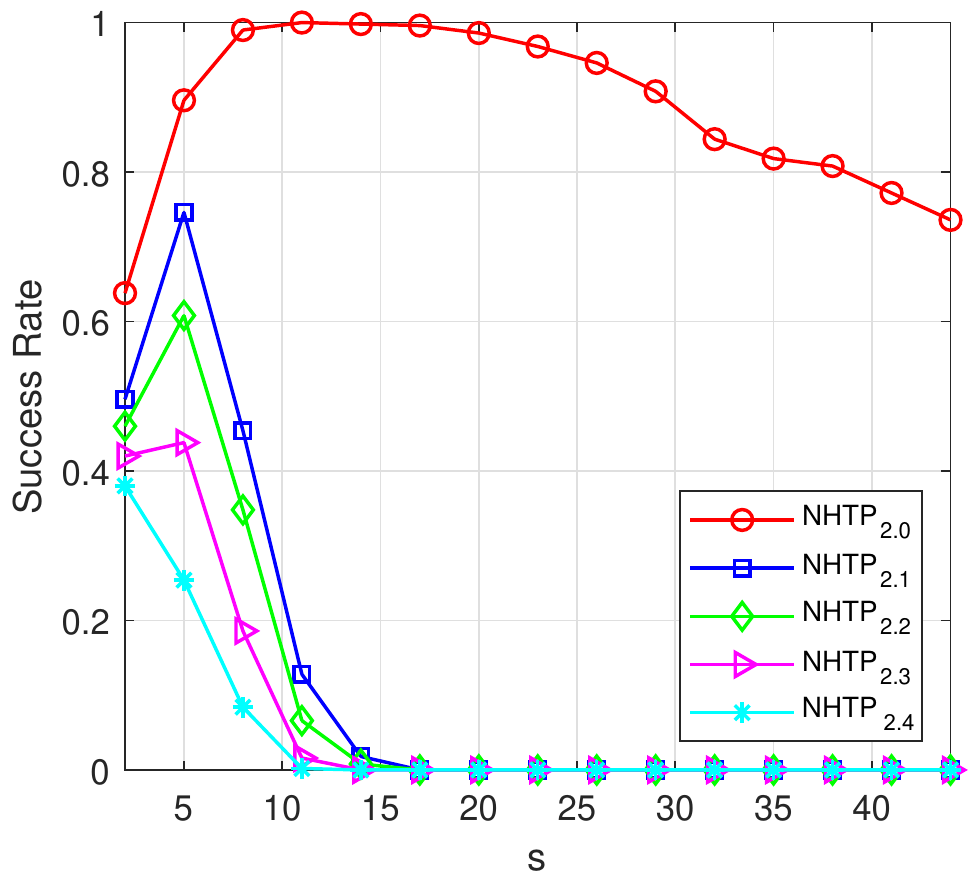}
  \caption{  \Cref{non-sdp-matrix}}
\end{subfigure}\vspace{-5mm}
\caption{Success rates of \NHTP$_r$. $n = 200,  s\in\{2,5,\cdots, 44\}$.}
\label{fig:SuccRate-r}
\end{figure}

Results for  \Cref{sdp-matrix} are presented in  \Cref{fig:SuccRate-r} (a), where $_r$ is set as $r=2,2.5,3,3.5,4$. It can be clearly seen that success rates decrease  along with $r$ ascending. We also test  other choices of $r=2.1,2.2,2.3,2.4$ and their results are between the red and blue lines with similar declined trends. For  \Cref{non-sdp-matrix}, we show success rates in  \Cref{fig:SuccRate-r} (b) generated by \NHTP$_r$ with $r=2,2.1,2.2,2.3,2.4$. We also tested  \NHTP$_r$ with $r>2.4$ and corresponding success rates are smaller than the case of $r=2.4$. Again, \NHTP$_{2.0}$ performs much better than the others. For each $s=s^*$, success rates decrease when $r$ ascends.  In conclusion, for fixed $s$, the smaller $r$ is (or for fixed $r$, the smaller $s$ is), the better recovery ability of \NHTP$_r$ has. 
\vspace{-5mm}

  \begin{table}[H]
  \centering
  \caption{Comparison of \NHTP$_r$  with different $r$.}\label{nhtp-r-z}
{\tiny\begin{tabular}{lcccccccccccc}
  \hline
  &&\multicolumn{5}{c}{$\|x-x^*\|/\|x^*\|$}&&\multicolumn{5}{c}{Time (seconds)} \\\cline{3-7} \cline{9-13}
  \multicolumn{13}{c}{ \Cref{z-matrix}}\\\hline
$n$ &&	5000	&	10000	&	15000	&	20000 &	25000	&&	5000	&	10000	&	15000	&	20000 &	25000	\\\hline
\NHTP$_{2.0}$	&&	0.00e-0	&	0.00e-0	&	0.00e-0	&	0.00e-0	&	0.00e-0	&&	0.004	&	0.005	&	0.007	&	0.009	&	0.014	\\
\NHTP$_{2.5}$	&&	5.65e-6	&	5.65e-6	&	5.65e-6	&	5.65e-6	&	5.65e-6	&&	0.011	&	0.031	&	0.021	&	0.027	&	0.038	\\
\NHTP$_{3.0}$	&&	2.44e-4	&	2.44e-4	&	2.44e-4	&	2.44e-4	&	2.44e-4	&&	0.011	&	0.021	&	0.025	&	0.029	&	0.038	\\
\NHTP$_{3.5}$	&&	7.84e-4	&	7.84e-4	&	7.84e-4	&	7.84e-4	&	7.84e-4	&&	0.014	&	0.022	&	0.028	&	0.039	&	0.045	\\
\NHTP$_{4.0}$	&&	2.28e-3	&	2.28e-3	&	2.28e-3	&	2.28e-3	&	2.28e-3	&&	0.018	&	0.062	&	0.031	&	0.040	&	0.047	\\\hline
 \multicolumn{13}{c}{ \Cref{sdp-matrix}}\\\hline
\NHTP$_{2.0}$	&	&	5.8e-12	&	6.4e-10	&	1.1e-12	&	1.4e-11	&	1.1e-10	&	&	0.064	&	0.188	&	0.394	&	0.699	&	1.021	\\
\NHTP$_{2.5}$	&	&	3.61e-5	&	1.63e-6	&	9.88e-7	&	4.03e-6	&	3.07e-5	&	&	0.113	&	0.381	&	0.788	&	1.567	&	2.431	\\
\NHTP$_{3.0}$	&	&	1.17e-2	&	9.03e-3	&	1.15e-2	&	4.12e-3	&	5.49e-3	&	&	0.173	&	0.540	&	1.169	&	2.423	&	4.075	\\
\NHTP$_{3.5}$	&	&	3.10e-2	&	1.37e-2	&	1.57e-2	&	1.13e-2	&	7.82e-3	&	&	0.203	&	0.684	&	1.958	&	3.649	&	6.031	\\
\NHTP$_{4.0}$	&	&	4.47e-2	&	2.25e-2	&	2.43e-2	&	4.92e-2	&	3.18e-2	&	&	0.227	&	0.800	&	2.563	&	5.150	&	7.179	\\
  \hline
\end{tabular}}
\end{table}
To see the accuracy of the solutions and the speed of \NHTP$_r$, we now test on two examples with higher dimensions $n$. For  \Cref{z-matrix}, we increase $n$ from $5000$ to $25000$ and fix $s^*=1$.   Results are presented in   \Cref{nhtp-r-z}.  Whilst for \Cref{sdp-matrix}, we
run independent 20 trials for each  $(n, s^*)$ with $n$ ranging  from $2000$ to $10000$ and  fixing $s^*=0.01n$. Average results {{over 20 trials}} are presented in  \Cref{nhtp-r-z}.  Clearly, for both examples, \NHTP$_{2}$ gets the most accurate solutions and runs the fastest for all cases. In a nutshell, the smaller $r$ is, the better \NHTP\ performs.

 \subsection{Effect of $s$ with fixing $r=2$}\label{sub:effect-s}
\begin{figure} [H]
\centering
\begin{subfigure}{0.49\textwidth}\centering
  \includegraphics[width=.9\linewidth]{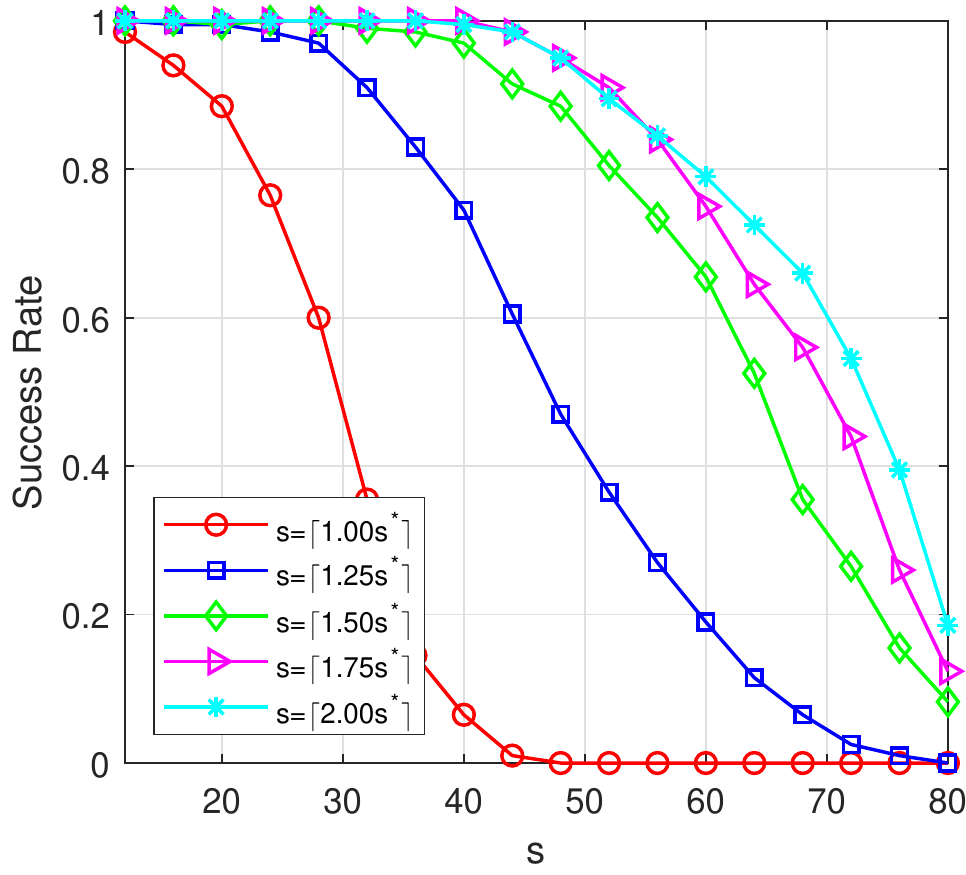}
   \caption{ \Cref{sdp-matrix}}
\end{subfigure}
\begin{subfigure}{0.49\textwidth}\centering
  \includegraphics[width=.9\linewidth]{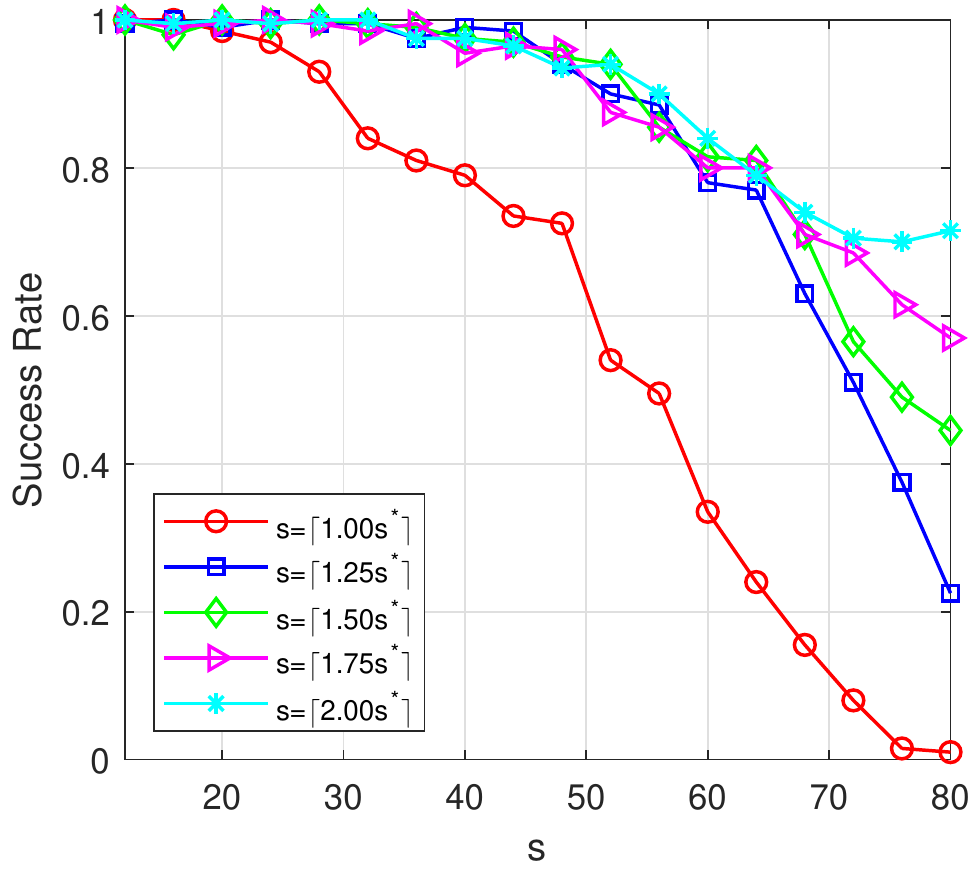}
  \caption{ \Cref{non-sdp-matrix}}
\end{subfigure}\vspace{-5mm}
\caption{Success rates of \NHTP$_2$. $n = 200,  s\in\{12,16,\cdots, 80\}$.}
\label{fig:SuccRate-s}\vspace{-5mm}
\end{figure}

 To make results comparable, we fix $r=2$.  In $S$, there is a parameter $s$ that should be given in advance. However, it is difficult to set an exact value for $s$ in practice. To see how the choices of {{$s$ affect the solution}} to (\ref{sncp}), we apply \NHTP\  to address  three examples with different$$ {{s\in \{s^*,\lceil 1.25s^*\rceil,\lceil 1.5s^*\rceil,\lceil 1.75s^*\rceil,2s^*\},}}$$ where $\lceil a\rceil$ returns the smallest integer that is no less than $a$.  To see the recovery ability, we first apply them  to solve  \Cref{sdp-matrix} and \Cref{non-sdp-matrix} with fixing $n = 200$  but with increasing sparsity level $s^*$ from $12$ to $80$. For each $(n, s^*)$, we run $500$ independent trials and record the corresponding success rates in   \Cref{fig:SuccRate-s}, where data show  that \NHTP$_2$ with $s>s^*$  generates better success rates than $s=s^*$. More detailed, the larger $s$ is, the higher success rates are produced by \NHTP$_2$. In addition, {{it seems to be more difficult for \NHTP$_r$ to solve \Cref{sdp-matrix} than \Cref{non-sdp-matrix}.}}  {{For instance}}, when $s=40$, \NHTP$_2$ is able to recover $80\%$ trials for  \Cref{non-sdp-matrix} while only get $5\%$  trials for \Cref{sdp-matrix}.

\begin{table}[H]
  \centering
  \caption{Comparison of \NHTP$_2$  with different $s$.}\label{nhtp-r-z-s}
{\tiny\begin{tabular}{lccccccccccc}
  \hline
  &\multicolumn{5}{c}{$\|x-x^*\|/\|x^*\|$}&&\multicolumn{5}{c}{Time (seconds)} \\\cline{2-6} \cline{8-12}
    \multicolumn{12}{c}{ \Cref{z-matrix}}\\\hline
$s\setminus n$ &	5000	&	10000	&	15000	&	20000 &	25000	&&	5000	&	10000	&	15000	&	20000 &	25000	\\\hline
$\lceil 1.00s^*\rceil$		&	0.0e-0	&	0.0e-0	&	0.0e-0	&	0.0e-0	&	0.0e-0	&	&	0.007	&	0.009	&	0.010	&	0.012	&	0.012	\\
$\lceil 1.25s^*\rceil$		&	0.0e-0	&	1.1e-16	&	0.0e-0	&	3.4e-21	&	0.0e-0	&	&	0.009	&	0.013	&	0.013	&	0.016	&	0.016	\\
$\lceil 1.50s^*\rceil$		&	0.0e-0	&	1.1e-16	&	0.0e-0	&	3.4e-21	&	0.0e-0	&	&	0.008	&	0.014	&	0.013	&	0.016	&	0.016	\\
$\lceil 1.75s^*\rceil$		&	0.0e-0	&	1.1e-16	&	0.0e-0	&	3.4e-21	&	0.0e-0	&	&	0.009	&	0.013	&	0.014	&	0.015	&	0.016	\\
$\lceil 2.00s^*\rceil$		&	0.0e-0	&	1.1e-16	&	0.0e-0	&	3.4e-21	&	0.0e-0	&	&	0.008	&	0.015	&	0.014	&	0.016	&	0.017	\\
  \hline
    \multicolumn{12}{c}{ \Cref{sdp-matrix}}\\\hline
$\lceil 1.00s^*\rceil$		&	6.6e-13	&	4.8e-11	&	5.5e-12	&	7.3e-13	&	9.2e-11	&	&	0.05 	&	0.19 	&	0.38 	&	0.68 	&	1.10 	\\
$\lceil 1.25s^*\rceil$		&	5.4e-13	&	1.6e-10	&	1.1e-16	&	1.0e-13	&	3.9e-15	&	&	0.06 	&	0.19 	&	0.41 	&	0.70 	&	1.10 	\\
$\lceil 1.50s^*\rceil$		&	1.5e-14	&	9.2e-11	&	1.8e-14	&	5.4e-14	&	8.7e-16	&	&	0.06 	&	0.19 	&	0.42 	&	0.71 	&	1.13 	\\
$\lceil 1.75s^*\rceil$		&	4.6e-11	&	1.2e-10	&	2.3e-16	&	6.5e-14	&	9.3e-16	&	&	0.06 	&	0.21 	&	0.43 	&	0.74 	&	1.16 	\\
$\lceil 2.00s^*\rceil$		&	4.7e-11	&	1.2e-10	&	3.4e-14	&	2.2e-14	&	2.1e-15	&	&	0.06 	&	0.21 	&	0.45 	&	0.77 	&	1.21 	\\\hline
\end{tabular}}
\end{table}

We now increase $n$ from $5000$ to $25000$ and fix $s^*=1$ for  \Cref{z-matrix}. Related results are presented in  \Cref{nhtp-r-z-s}.  While for  \Cref{sdp-matrix}, we again run independent 20 trials for each  $(n, s^*)$ with $n$ ranging  from $2000$ to $10000$ and  keeping $s^*=0.01n$. Average results {{over 20 trials}} are presented in \Cref{nhtp-r-z-s}. For both tables, it can be clearly seen that accuracies obtained by \NHTP$_2$ under different $s$ are similar. {{As  expected, smaller $s$ enables \NHTP$_2$ to run slightly faster than larger $s$.}}

 \subsection{Strategies to select $s$} {{Assume the sparse LCP \eqref{sncp} admits a sparsest solution $x^*$ with sparsity level $s^*$. As long as $s^*\ll n$ (e.g. $s^*\leq \lceil0.1n\rceil$), numerical experiments in \Cref{sub:effect-s} demonstrate  that \NHTP\ achieves the sparsest solutions with a very high possibility if we set $s\geq s^*$, see \Cref{fig:SuccRate-s} for instance. Therefore, a possible way to tune a proper $s$ is designed as \Cref{s-tuning}, where parameter can be set as $s_0=\lceil  n/5000 \rceil$,  $\varrho=\max\{2,\log_{10}(n)\}$ and $\epsilon=10^{-8}$. In this way, if \eqref{sncp} admits a solution $x^*$ with $s^*\ll n$, then the worst case to achieve $s\geq s^*$ is running \NHTP\ $\lceil \log_{\varrho}(s^*/s_0)\rceil$ times, after which \NHTP\ will possibly achieve  the solution. 

\begin{algorithm}
\caption{{\tt NHTPT}: \NHTP\ with sparsity level tuning}
\label{s-tuning}
\begin{algorithmic}
\STATE{Initialize a small integer $s_0\in\N, \varrho>1, \epsilon>0$ and $x^0=0$. Set $\ell\Leftarrow0$.}
\WHILE{$f_r(x^\ell)\geq \epsilon$}
\STATE{ Set $s=s_\ell$ and run \NHTP\ in \Cref{alg:nhtp} to generate a solution $x^{\ell+1}$}
\STATE{ Set $s_\ell=\lceil  \varrho s_\ell\rceil$ and $\ell\Leftarrow\ell+1.$}
\ENDWHILE
\RETURN the solution $x^{\ell}$.
\end{algorithmic}
\end{algorithm}

An alternative takes advantage of other methods that do not need the prior information $s$, for example,  Lemke's ({\tt Lemke} \footnote{available at \url{http://ftp.cs.wisc.edu/math-prog/matlab/lemke.m}}) algorithm, a well-known  high standard method to solve the LCP. Therefore, we could first run {\tt Lemke} to obtain a solution $x_{\rm lem}$ and then set $s=\|x_{\rm lem}\|_0$ for \NHTP. Note that $\|x_{\rm lem}\|_0$ actually provides an upper bound of $s$. However, we test  that this upper bound sometimes is good enough.

Now we would like to see the performance of {\tt Lemke}, \NHTP\ with the help of  $s=\|x_{\rm lem}\|_0$ and {\tt NHTPT} in \Cref{s-tuning}. We fix $r=2$ in $f_r$ for the latter two methods.  Average results over 20 trials are presented in \Cref{nhtp-lemke}, where all methods achieve solutions to LCP for all cases since the objective values $f_r$ are close to zeros. 
For \Cref{sdp-matrix}, where the `ground truth' solutions are given and $s^*$ is set as $\lceil0.01n\rceil$, three methods render  solutions with sparsity levels being identical to $s^*$. \NHTP\ runs the fastest, followed by  {\tt NHTPT}. While {\tt Lemke} consumes too much time, e.g., 78.27 seconds v.s. 7.3 seconds by \NHTP\ when $n=25000$. For \Cref{non-sdp-matrix-nox},  the `ground truth' solutions are unknown and $s^*$ is set as $\lceil 0.5n\rceil$. Note that this large $s^*$ for such example is not the sparsity level of a solution, but can be an upper bound of $s$. As shown in \Cref{nhtp-lemke},  three methods succeed in finding very sparse solutions since the sparsity levels $\|x\|_0$ are relatively small to the large $s^*$.  In addition, {\tt NHTPT} runs the fastest and also produces the sparsest solutions, followed by \NHTP. 

The performance of  {\tt NHTPT} solving the above two examples  illustrates that the strategy in \Cref{s-tuning} allows \NHTP\ to find a proper $s$ iteratively. However, in the sequel, we still focus on \NHTP\ itself instead of {\tt NHTPT} for the sake of simplicity.
 
\begin{table} 
  \centering
  \caption{Comparison of {\tt Lemke}, \NHTP\ and {\tt NHTPT}.}\label{nhtp-lemke}
{\tiny\begin{tabular}{lcccccccccccc}
  \hline
  &\multicolumn{3}{c}{$f_2$}&&\multicolumn{3}{c}{Time (seconds)}&&\multicolumn{3}{c}{$\|x\|_0$}\\
$n$	&	{\tt Lemke}	&	\NHTP	&	{\tt NHTPT}	&	&	{\tt Lemke}	&	\NHTP	&	{\tt NHTPT}	&	&	{\tt Lemke}	&	\NHTP	&	{\tt NHTPT}	\\\cline{2-4}\cline{6-8}  \cline{10-12}    
 &\multicolumn{11}{c}{\Cref{sdp-matrix} } \\ \hline
5000	&	6.63e-30	&	5.22e-15	&	3.78e-14	&	&	0.63 	&	0.27 	&	0.61 	&	&	50 	&	50 	&	50 	\\
10000	&	1.32e-29	&	2.25e-14	&	9.35e-15	&	&	3.81 	&	0.95 	&	2.05 	&	&	100 	&	100 	&	100 	\\
15000	&	3.09e-29	&	7.36e-15	&	3.68e-15	&	&	12.1 	&	2.03 	&	4.41 	&	&	150 	&	150 	&	150 	\\
20000	&	4.93e-29	&	6.77e-15	&	5.39e-15	&	&	27.6 	&	3.55 	&	7.90 	&	&	200 	&	200 	&	200 	\\
25000	&	1.09e-28	&	4.82e-14	&	3.23e-16	&	&	78.3 	&	7.30 	&	12.7  	&	&	250 	&	250 	&	250 	\\\hline
   &\multicolumn{11}{c}{\Cref{non-sdp-matrix-nox} } \\ \hline
5000	&	3.63e-09	&	1.50e-12	&	2.45e-11	&	&	0.43 	&	0.28 	&	0.16 	&	&	25.7 	&	25.7 	&	1.0 	\\
10000	&	1.23e-08	&	8.35e-11	&	7.25e-12	&	&	1.29 	&	0.62 	&	0.48 	&	&	21.4 	&	21.4 	&	2.0 	\\
15000	&	4.44e-09	&	5.44e-12	&	2.64e-12	&	&	2.90 	&	1.16 	&	1.10 	&	&	10.8 	&	10.8 	&	2.9 	\\
20000	&	9.87e-09	&	1.12e-12	&	2.04e-12	&	&	5.21 	&	1.88 	&	1.79 	&	&	6.3 	&	6.3 	&	4.0 	\\
25000	&	1.65e-08	&	2.14e-12	&	1.20e-12	&	&	30.2 	&	4.90 	&	2.76 	&	&	5.9 	&	5.9 	&	4.9 	\\

\hline
\end{tabular}}
\end{table} 

}}

 \subsection{Numerical comparisons}
 Since there are very few methods that have been proposed to process the sparse LCP, we compare \NHTP$_r$ only with half thresholding projection (\HTP) method \cite{SZX15} and extra-gradient thresholding algorithm (\ETA) \cite{SZPX15}. We use all their default parameters and terminate both of them when $\|x^k-z^k\|<10^{-5}\max\{1,\|x^k\|\} $ or the maximum number of iterations reach 2000. Note that both methods make use of the first order information of the involved functions and thus belong to {{the class of the first order methods}}.  \NHTP\ uses the origin as its default starting point. However, as a second order method, it is suggested to start from {{a local area around a solution}}. Therefore, we take advantage of the solution obtained by \HTP\ as the starting point of \NHTP. Under such circumstance, write  \NHTP$_r$ as \HNHTP$_r$. We thus compare \NHTP$_2$, \HNHTP$_2$, \HNHTP$_{2.5}$, \HNHTP$_3$, \HTP\ and \ETA. For the former four \NHTP-related methods, we choose $s=s^*$ in $S$ for \Cref{z-matrix}, \Cref{sdp-matrix} and \Cref{non-sdp-matrix} since the sparsity of the `ground truth' solution is $s^*$ and choose
 $$s=\min\{\|x_{\HTP}\|_0,\|x_{\ETA}\|_0,s^*\}.$$
 for  \Cref{non-sdp-matrix-nox} since the `ground truth' solution is unknown, where $x_{\HTP}$ and $x_{\ETA}$ are solutions produced by \HTP\ and \ETA, respectively. In such a way, \NHTP\ {{could always get}} solutions that are sparser than  solutions produced by the last two methods.

  \begin{figure}[H]
\centering
\begin{subfigure}{0.49\textwidth}\centering
  \includegraphics[width=.9\linewidth]{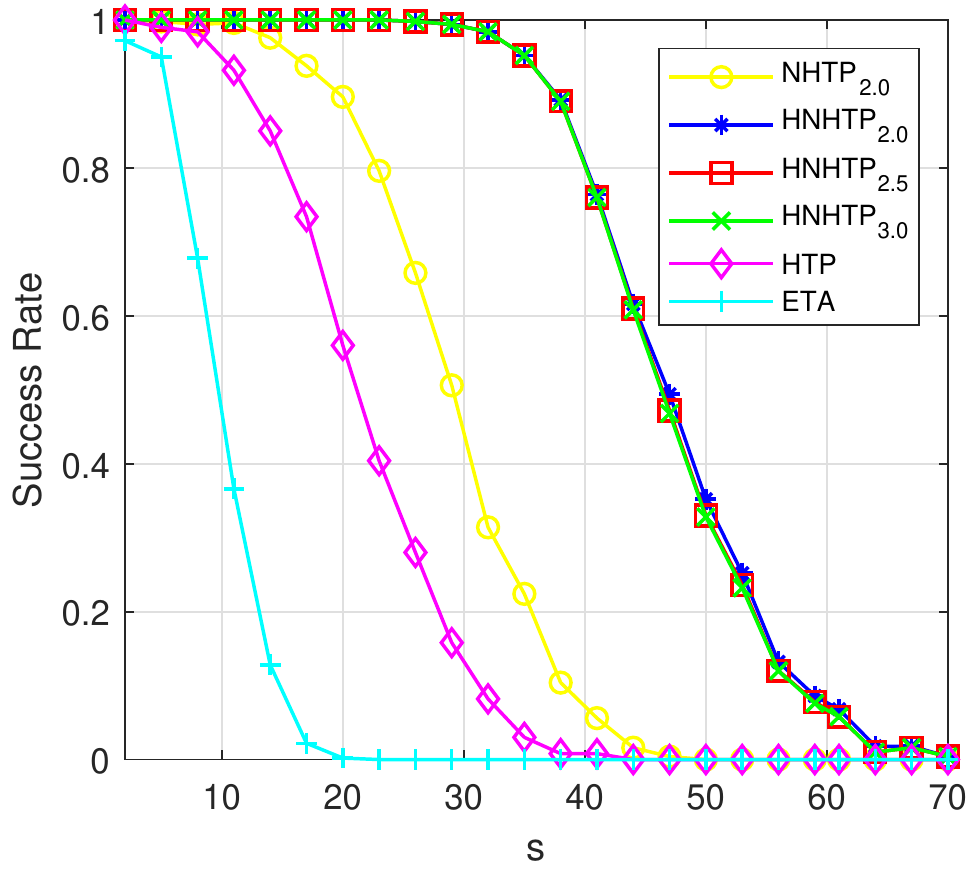}
   \caption{ \Cref{sdp-matrix}}
\end{subfigure}
\begin{subfigure}{0.49\textwidth}\centering
  \includegraphics[width=.9\linewidth]{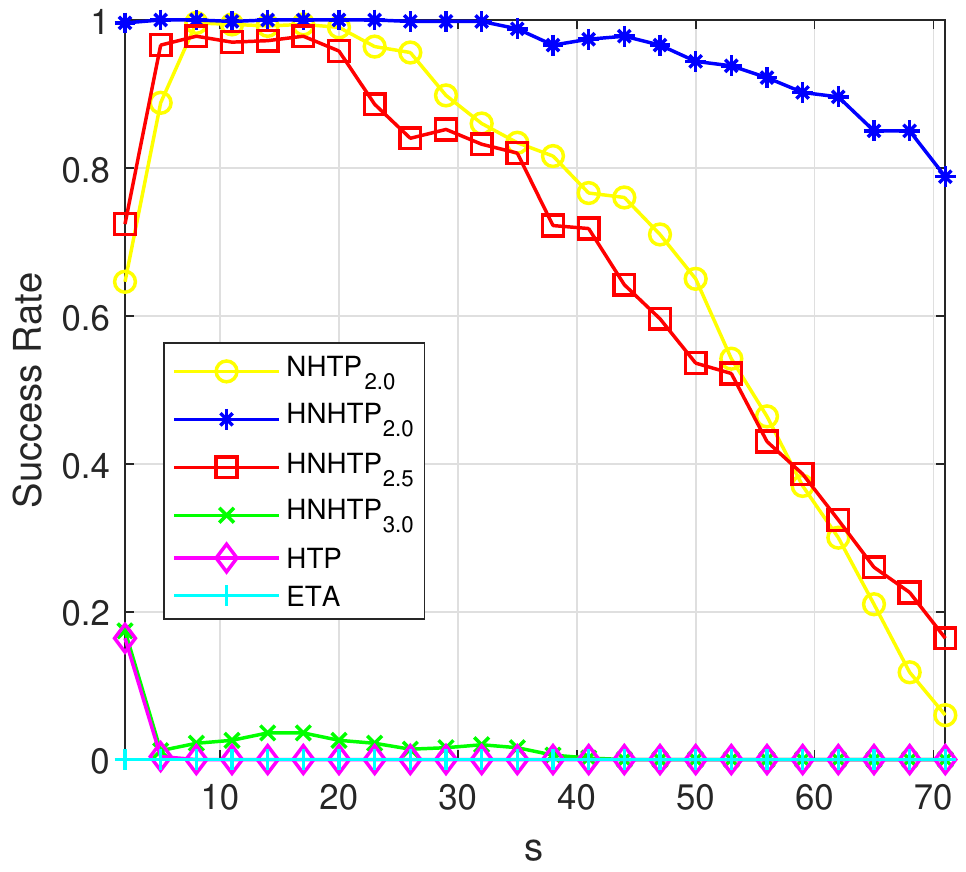}
  \caption{ \Cref{non-sdp-matrix}}
\end{subfigure}\vspace{-5mm}
\caption{Success rates of \NHTP, \HTP\ and \ETA. $n = 200, s\in\{2,5,\cdots,71\}$.}
\label{fig:SuccRate-6}
\end{figure}

{\bf (a) Recovery ability.} Similarly, to see the recovery ability, we first apply them  to solve  \Cref{sdp-matrix} and \Cref{non-sdp-matrix} with fixing $n = 200$  but with increasing sparsity level $s^*$ from $2$ to $71$. For each $(n, s^*)$, we run $500$ independent trials and record the corresponding success rates in  \Cref{fig:SuccRate-6},  where data in subfigure (a)  show that \HNHTP$_2$, \HNHTP$_{2.5}$, \HNHTP$_3$ generate similar results and obtain the highest success rates, followed by \NHTP$_2$. While \HTP\ and \ETA\ come the last. When those methods are applied to solve \Cref{non-sdp-matrix}, the results in subfigure (b) present a big different picture.   \HNHTP$_2$  outperforms the other five methods, followed by \NHTP$_2$, \HNHTP$_{2.5}$. {{In contrast}}, \HNHTP$_3$ \HTP\ and \ETA\ basically fail to recover solutions for  cases of $s\geq5.$   Overall, one could conclude that \HTP\ itself does not produce accurate solutions but could offer good starting points,  from which \HNHTP$_2$, \HNHTP$_{2.5}$, \HNHTP$_3$ benefit significantly.

{\bf (b) Accuracy and speed in the higher dimensional setting.} To see the performance of six methods on solving larger size problems, we now increase $n$ from $5000$ to $25000$ and fix $s^*=1$ for  \Cref{z-matrix}. Related results are presented in  \Cref{nhtp-r-z-1}.  {{For  \Cref{sdp-matrix}, we again run independent 20 trials}} for each  $(n, s^*)$ with $n$ ranging  from $2000$ to $10000$ and  keeping $s^*=0.01n$. Average results {{over 20 trials}} are presented in  \Cref{nhtp-r-z-1}. It can be clearly seen that \HNHTP$_2$ and \NHTP$_2$ get the most accurate solutions, followed by \HNHTP$_{2.5}$ and \HNHTP$_3$, \HTP\ comes the last.  For the computational time, all \NHTP\ methods run much faster than \HTP\ and \ETA.


\begin{table}[H]
  \centering
  \caption{Comparison of \NHTP$_r$, \HTP\ and \ETA.}\label{nhtp-r-z-1}
{\tiny\begin{tabular}{lcccccccccccc} \hline
  &&\multicolumn{5}{c}{$\|x-x^*\|/\|x^*\|$}&&\multicolumn{5}{c}{Time (seconds)} \\\cline{3-7} \cline{9-13}
    \multicolumn{13}{c}{ \Cref{z-matrix}}\\\hline
$n$ &&	5000	&	10000	&	15000	&	20000 &	25000	&&	5000	&	10000	&	15000	&	20000 &	25000	\\\hline
\NHTP$_{2.0}$	&	&	0.00e-0	&	0.00e-0	&	0.00e-0	&	0.00e-0	&	0.00e-0	&	&	0.004	&	0.006	&	0.008	&	0.007	&	0.009	\\
\HNHTP$_{2.0}$	&	&	0.00e-0	&	0.00e-0	&	0.00e-0	&	0.00e-0	&	0.00e-0	&	&	0.003	&	0.005	&	0.007	&	0.006	&	0.009	\\
\HNHTP$_{2.5}$	&	&	3.89e-6	&	3.89e-6	&	3.89e-6	&	3.89e-6	&	3.89e-6	&	&	0.006	&	0.010	&	0.014	&	0.020	&	0.017	\\
\HNHTP$_{3.0}$	&	&	7.88e-5	&	7.87e-5	&	7.87e-5	&	7.87e-5	&	7.87e-5	&	&	0.004	&	0.006	&	0.007	&	0.012	&	0.010	\\
\HTP	&	&	3.15e-4	&	3.15e-4	&	3.15e-4	&	3.15e-4	&	3.15e-4	&	&	0.037	&	0.086	&	0.132	&	0.163	&	0.171	\\
\ETA	&	&	2.93e-4	&	2.93e-4	&	2.93e-4	&	2.93e-4	&	2.93e-4	&	&	0.077	&	0.193	&	0.282	&	0.498	&	0.378	\\
  \hline
  \multicolumn{13}{c}{ \Cref{sdp-matrix}}\\\hline
\NHTP$_{2.0}$	&	&	2.0e-12	&	2.5e-11	&	1.3e-12	&	1.0e-13	&	1.5e-13	&	&	0.08 	&	0.21 	&	0.40 	&	0.73 	&	1.02 	\\
\HNHTP$_{2.0}$	&	&	6.7e-12	&	1.3e-10	&	1.7e-11	&	1.0e-11	&	4.2e-11	&	&	0.05 	&	0.18 	&	0.35 	&	0.62 	&	0.98 	\\
\HNHTP$_{2.5}$	&	&	1.76e-7	&	8.51e-8	&	7.98e-8	&	7.95e-8	&	1.41e-7	&	&	0.10 	&	0.33 	&	0.66 	&	1.18 	&	1.90 	\\
\HNHTP$_{3.0}$	&	&	4.41e-6	&	1.61e-6	&	1.51e-6	&	2.71e-6	&	4.48e-6	&	&	0.09 	&	0.30 	&	0.60 	&	1.04 	&	1.70 	\\
\HTP	&	&	2.38e-4	&	3.17e-4	&	2.94e-4	&	2.85e-4	&	4.21e-4	&	&	1.61 	&	6.96 	&	13.7 	&	24.6 	&	41.8  	\\
\ETA	&	&	2.01e-4	&	1.76e-4	&	1.73e-4	&	1.69e-4	&	2.71e-4	&	&	3.07 	&	15.1 	&	30.9  	&	56.6 	&	88.4 	\\
  \hline
\end{tabular}}
\end{table}

 \begin{table} [H]
  \centering
  \caption{Comparison of \NHTP$_r$, \HTP\ and \ETA.}\label{nhtp-ex4}
{\tiny\begin{tabular}{lcccccccccccc}
  \hline
  &&\multicolumn{5}{c}{$f_2(x)$}&&\multicolumn{5}{c}{Time (seconds)} \\\cline{3-7} \cline{9-13}
$n$ &&	2000	&	4000	&	6000	&	8000 &	10000	&&	2000	&	4000	&	6000	&	8000 &	10000	\\\hline
\NHTP$_{2.0}$	&	&	9.65e-8	&	1.25e-7	&	3.08e-8	&	1.22e-7	&	1.10e-7	&	&	0.04	&	0.04	&	0.07	&	0.12	&	0.16	\\
\HNHTP$_{2.0}$	&	&	7.63e-6	&	1.79e-8	&	1.41e-8	&	1.29e-7	&	4.52e-8	&	&	0.07	&	0.04	&	0.05	&	0.07	&	0.11	\\
\HNHTP$_{2.5}$	&	&	6.08e-5	&	2.66e-5	&	1.13e-8	&	1.69e-8	&	1.63e-8	&	&	0.02	&	0.04	&	0.06	&	0.09	&	0.12	\\
\HNHTP$_{3.0}$	&	&	8.72e-5	&	3.33e-4	&	5.34e-4	&	7.07e-4	&	7.65e-4	&	&	0.01	&	0.02	&	0.02	&	0.04	&	0.05	\\
\HTP	&	&	9.66e-5	&	1.89e-4	&	2.86e-4	&	3.68e-4	&	4.50e-4	&	&	0.03	&	0.03	&	0.04	&	0.04	&	0.05	\\
\ETA	&	&	1.84e-4	&	1.94e-4	&	2.39e-4	&	2.50e-4	&	3.09e-4	&	&	0.27	&	1.51	&	3.65	&	7.24	&	12.2	\\
  \hline
  &&\multicolumn{5}{c}{$\| \nabla f_2(x)\|$}&&\multicolumn{5}{c}{$\|x\|_0$} \\\cline{3-7} \cline{9-13}																 \NHTP$_{2.0}$	&	&	1.14e-4	&	2.05e-5	&	1.37e-5	&	3.30e-5	&	2.60e-5	&	&	9.2	&	15.0	&	19.3	&	24.3	&	30.5	\\
\HNHTP$_{2.0}$	&	&	1.93e-4	&	1.70e-5	&	2.39e-5	&	7.26e-5	&	3.67e-5	&	&	9.0	&	14.8	&	19.3	&	24.3	&	30.5	\\
\HNHTP$_{2.5}$	&	&	1.86e-3	&	3.02e-3	&	3.90e-3	&	5.86e-3	&	7.58e-3	&	&	9.2	&	15.0	&	19.3	&	24.3	&	30.5	\\
\HNHTP$_{3.0}$	&	&	4.85e-2	&	1.20e-2	&	7.24e-2	&	1.92e-1	&	2.74e-1	&	&	9.2	&	15.0	&	19.3	&	24.3	&	30.5	\\
\HTP	&	&	2.12e-3	&	4.77e-3	&	7.45e-3	&	1.04e-2	&	1.30e-2	&	&	11.3	&	27.1	&	44.8	&	64.8	&	81.7	\\
\ETA	&	&	3.30e-3	&	4.87e-3	&	6.68e-3	&	8.11e-3	&	1.02e-2	&	&	9.2	&	15.0	&	19.3	&	24.3	&	30.5	\\
 \hline
\end{tabular}}
\end{table}

{\bf (c) Performance on solving examples without known solutions.} Now we compare those methods on solving  \Cref{non-sdp-matrix-nox}, where solutions are unknown. Nevertheless,  they possibly admit some sparse solutions by \Cref{existence}. We run independent 20 trials for each  $(n, s^*)$ with $n$ ranging  from $2000$ to $10000$ and  keeping $s^*=0.01n$. Average results are presented in  \Cref{nhtp-ex4}. Note that since the objective functions $f_r$ is different with different $r$, to make comparison reasonable, we calculate $f_2(x )$, where $x$ is generated by one of six methods. For  \Cref{non-sdp-matrix-nox}, all \NHTP-related methods get the smallest objective function values and $ \| \nabla f_2(x)\|$ with the sparsest solutions, which means they outperform \HTP\ and \ETA\ in terms of the quality of solutions. In addition, \HTP\ always obtains solutions that are least sparse, but it and \HNHTP$_{3.0}$ run the fastest. \ETA\ is the slowest one again.

\subsection{Comparison of different NCP functions}\label{subsec:ncp}
For the sake of illustrating the advantage of $\phi_r$, we make use of \NHTP\ to address the problem (\ref{slco}) with different objective functions constructed by three NCP functions $\phi_{FB}$, $\phi_{\min}$ and $\phi^2_{II}$  from  \Cref{remark1}. The corresponding merit functions are
\begin{eqnarray*}
f_{FB}(x)&=&0.5\sum (\phi_{FB}(x_i,y_i))^2\\
 &=&0.5\left[\|x\|^2+\|y\|^2+\|x+y\|^2-2\langle \sqrt{x\circ x+y\circ y},  x+y\rangle\right],\\
f_{\min}(x)&=&0.5\sum (\phi_{\min}(x_i,y_i))^2\\
&=& 0.5 [ \|x+y\|^2+ \|x-y\|^2 - 2\langle  \sqrt{ (x-y)\circ (x-y) },  x+y\rangle ],\\
f_{II}(x)&=&0.5\sum (\phi^2_{II}(x_i,y_i))^2 = 0.5\left[\|(x\circ y)_+\|^2+\|x_-\|^2+\|y_-\|^2\right],
  \end{eqnarray*}
  where $\sqrt{z}=(\sqrt{z_1},\cdots,\sqrt{z_n})^\top $ and $y=Mx+q$.

  \begin{remark} We have some comments about the above merit functions and  $f_2$.
  \begin{itemize}
  \item[i)] Note that $f_{FB}$ and $f_{\min}$ have unbounded Hessian at $(0,0)$ and $x=y$, respectively. Therefore, to make use of \NHTP, we add a small scalar $\varepsilon$ (e.g. $10^{-10}$) to smooth $\sqrt{z}$, namely, replacing $\sqrt{z}$ by $\sqrt{z+\varepsilon}$ in $f_{FB}$ and $f_{\min}$. Then their gradients and Hessian are able to be derived. In addition, similar rules to calculate $\partial^2 f_2(x)$ in (\ref{hess-phi2-lcp}) also {{lead to}} the generalized Hessian of $f_{II}$.
\item[ii)] As shown in \cite{ZXQ19}, to derive the Newton direction, each step in \NHTP\  calculates a submatrix $\nabla_{TT}^2 f(x)$ of  the Hessian of $f$. It is easy to see that the Hessians of  $f_{FB}$ and $f_{\min}$ have a term $M^\top M$. Therefore, {{we need to }}compute $M_T^\top M_T$ and the computational complexity is about $\mathcal{O}(ns^2)$. While for $f_{II}$ and $f_{2}$, the most expensive computation is $ M_T^\top \Diag (\zeta )  M_T$. When the point {{is close to}} a solution to the LCP, then $y\geq0$, which together with (\ref{zetalcp}) indicates $$M_T^\top \Diag (\zeta )  M_T=M_{\supp(x)T}^\top \Diag (\zeta_{\supp(x)} )  M_{\supp(x)T}.$$
This means the computational complexity is about $\mathcal{O}(s^3)$. {{Therefore, we expect that $f_{FB}$ and $f_{\min}$ take longer time to do computations than $f_{II}$ and $f_{2}$ in each step, which is testified by the  numerical experiments in the sequel.}}
    \end{itemize}
  \end{remark}
 
\begin{figure} 
\begin{subfigure}{0.49\textwidth}
\begin{center}
  \includegraphics[width=.995\linewidth]{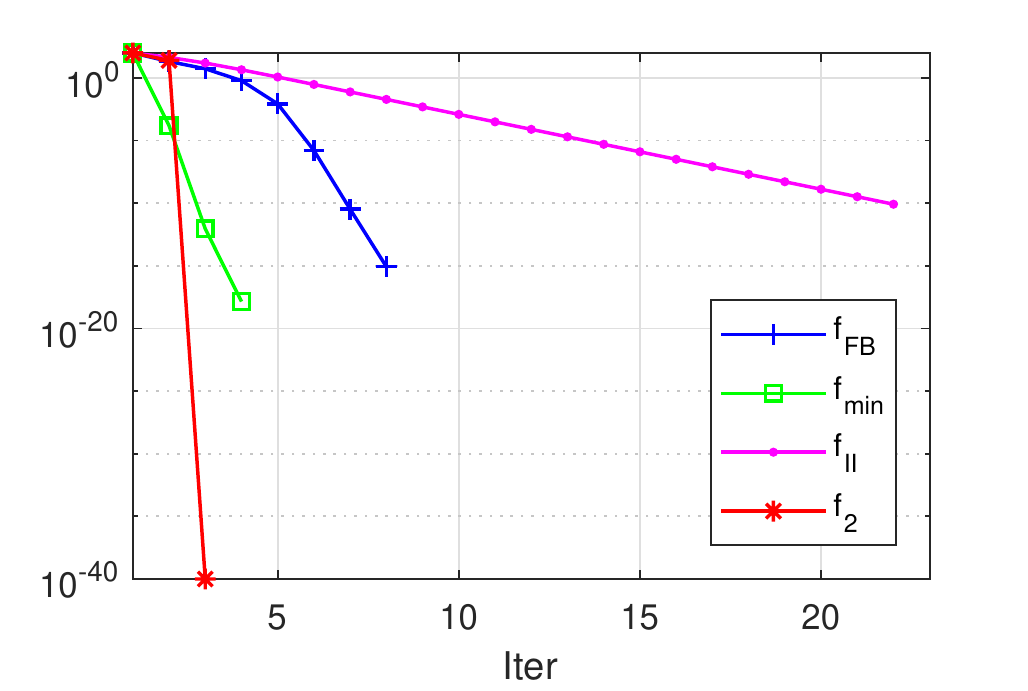}
  \caption{ \Cref{z-matrix}}\end{center}   
\end{subfigure}
\begin{subfigure}{0.49\textwidth}
 \begin{center} \includegraphics[width=.995\linewidth]{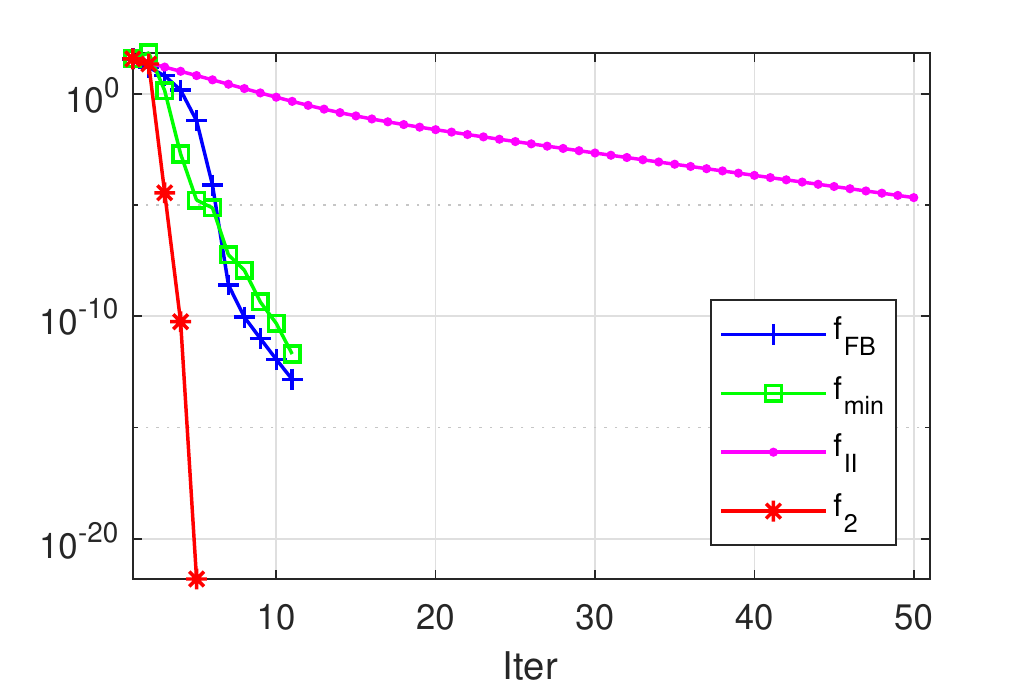}
  \caption{  \Cref{sdp-matrix}}\end{center} 
\end{subfigure} 
\begin{subfigure}{0.49\textwidth}
 \begin{center} \includegraphics[width=.995\linewidth]{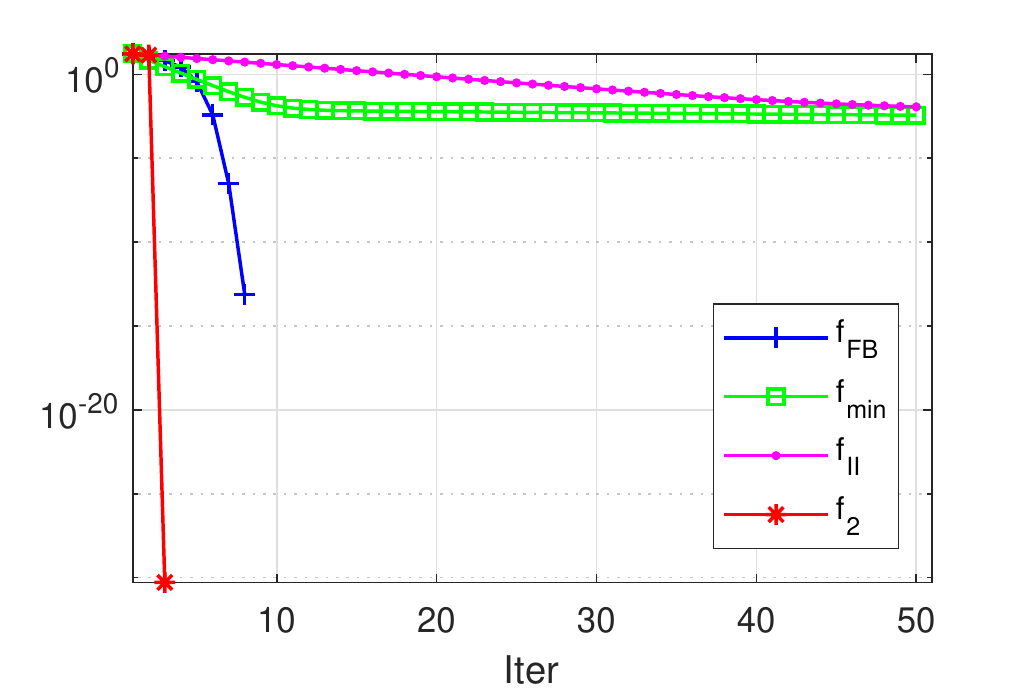}
   \caption{ \Cref{non-sdp-matrix}}\end{center} 
\end{subfigure}
\begin{subfigure}{0.49\textwidth}
 \begin{center} \includegraphics[width=.995\linewidth]{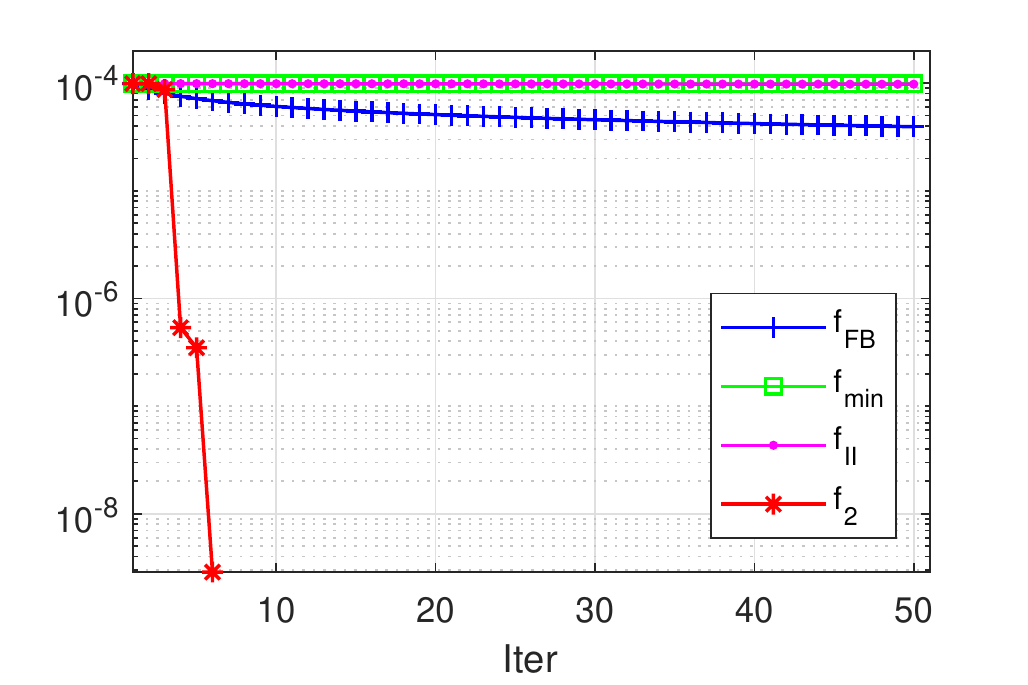}
  \caption{ \Cref{non-sdp-matrix-nox}}\end{center} 
\end{subfigure} \vspace{-5mm}
\caption{Objective function values $f_2$ at first 20 or 50 iterations. $n = 200,  s=2$.}
\label{fig:decreasing-obj}
\vspace{-5mm}
\end{figure}

 Now we apply \NHTP\ with fixing $s=s^*=0.01n$  to process the  sparsity constrained  model (\ref{slco})  with four merit functions $f_{FB}$, $f_{\min}$, $f_{II}$ and $f_2$. To see the decline of objective function values in each step at the beginning of the method, we report $f_2(x)$ to make results comparable, where $x$ is generated by \NHTP\ solving sparsity constrained  model with one of there merit  functions. For example,  we record the iterates $x^1,x^2,\cdots$ generated by \NHTP\ under $f_{FB}$ and then calculate $f_2(x^1), f_2(x^2),\cdots$. Results are presented in  \Cref{fig:decreasing-obj}. It is worth mentioning that all merit functions make \NHTP\ get the global solutions eventually, while we only report results at first 22 or 50 iterations. The prominent feature of the four sub-figures is that the lines of $f_2$ drop dramatically for all examples.  It only takes less than five steps to reduce the objective almost to zero. By contrast, when \NHTP\ addresses the model with $f_{II}$,  much more steps are required and the objective function values decline  relatively slowly. This phenomenon also appears for \Cref{non-sdp-matrix-nox}, where \NHTP\ seems not to prefer the sparsity constrained   models with $f_{FB}$, $f_{\min}$ and $f_{II}$.

    \begin{figure}[H] 
\centering
\begin{subfigure}{0.33\textwidth}
  \includegraphics[width=.995\linewidth]{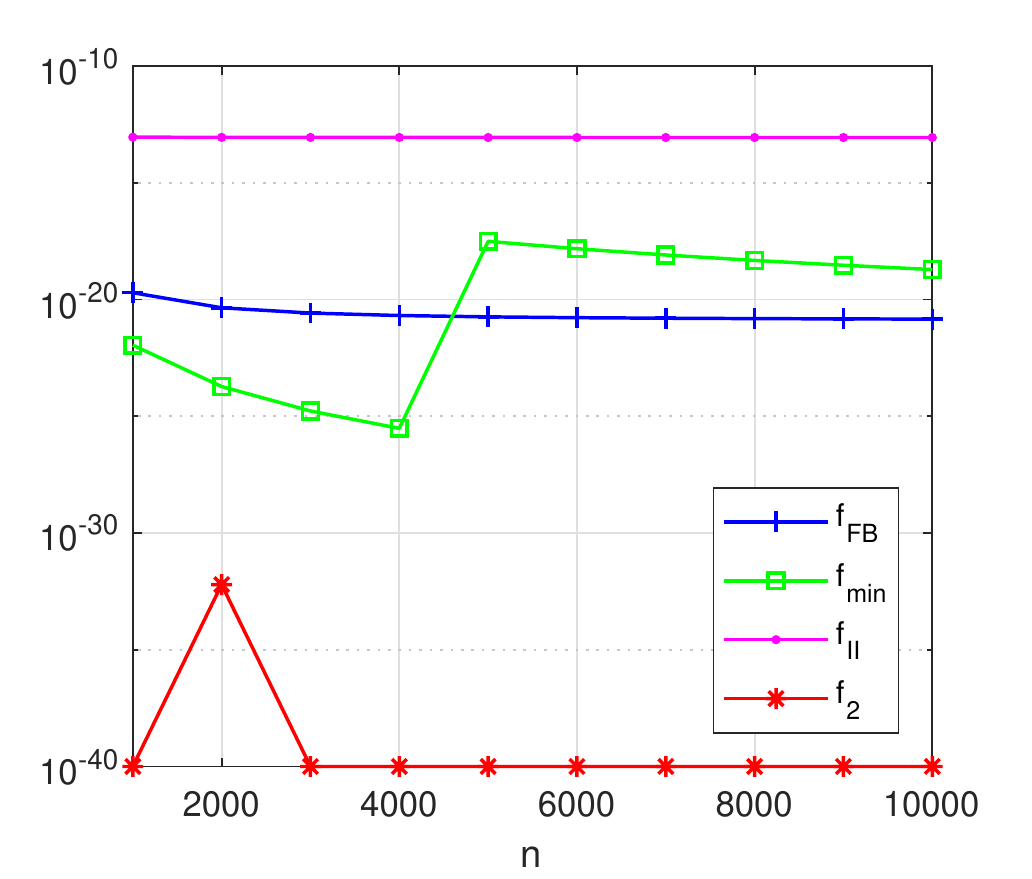}
   \caption{ \Cref{z-matrix}: $f_2(x)$}
\end{subfigure}
\begin{subfigure}{0.32\textwidth}
  \includegraphics[width=.995\linewidth]{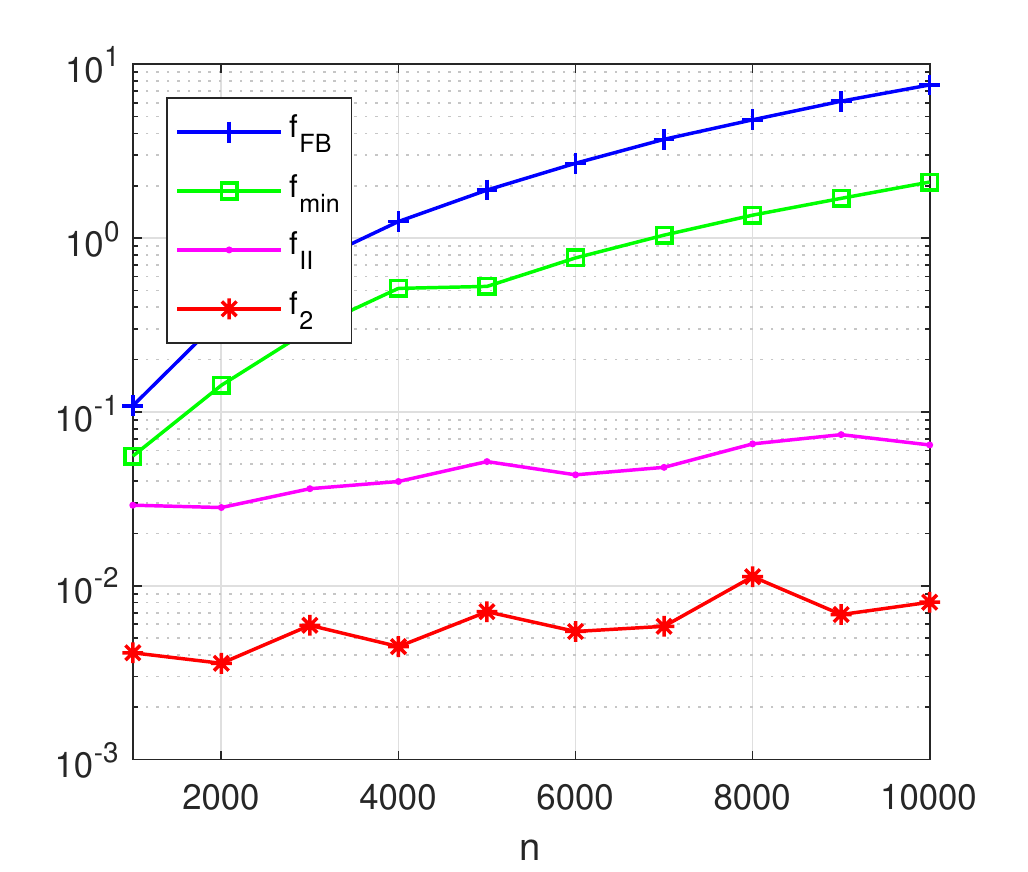}
  \caption{ \Cref{z-matrix}: Time}
\end{subfigure} 
\begin{subfigure}{0.32\textwidth}
  \includegraphics[width=.995\linewidth]{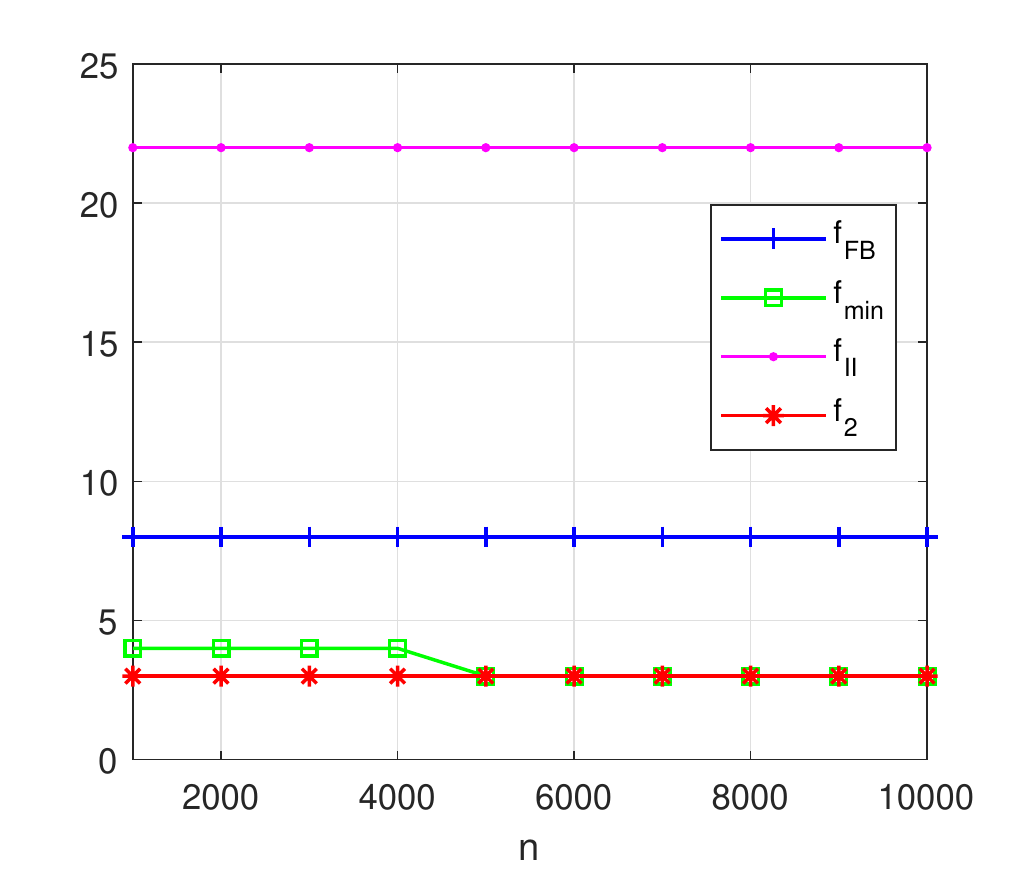}
   \caption{  \Cref{z-matrix}: Iteration}
\end{subfigure}\\
\begin{subfigure}{0.32\textwidth}
  \includegraphics[width=.995\linewidth]{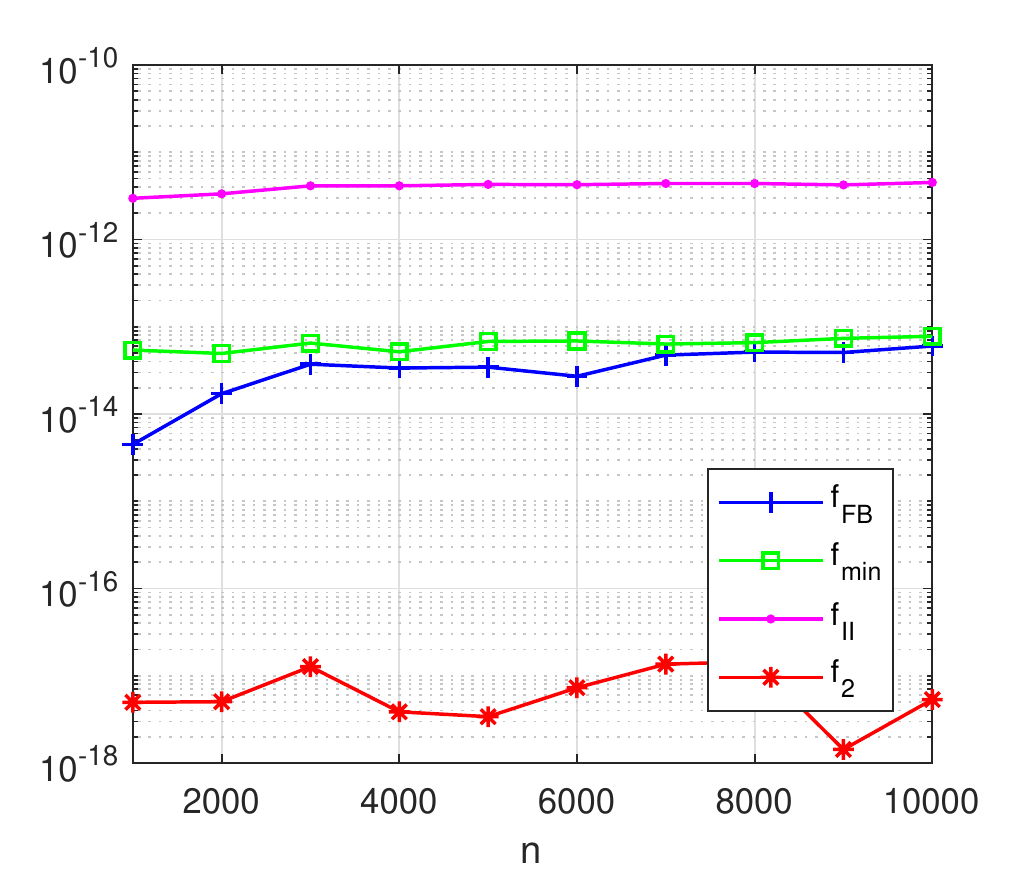}
  \caption{ \Cref{sdp-matrix}: $f_2(x)$}
\end{subfigure} 
\begin{subfigure}{0.32\textwidth}
  \includegraphics[width=.995\linewidth]{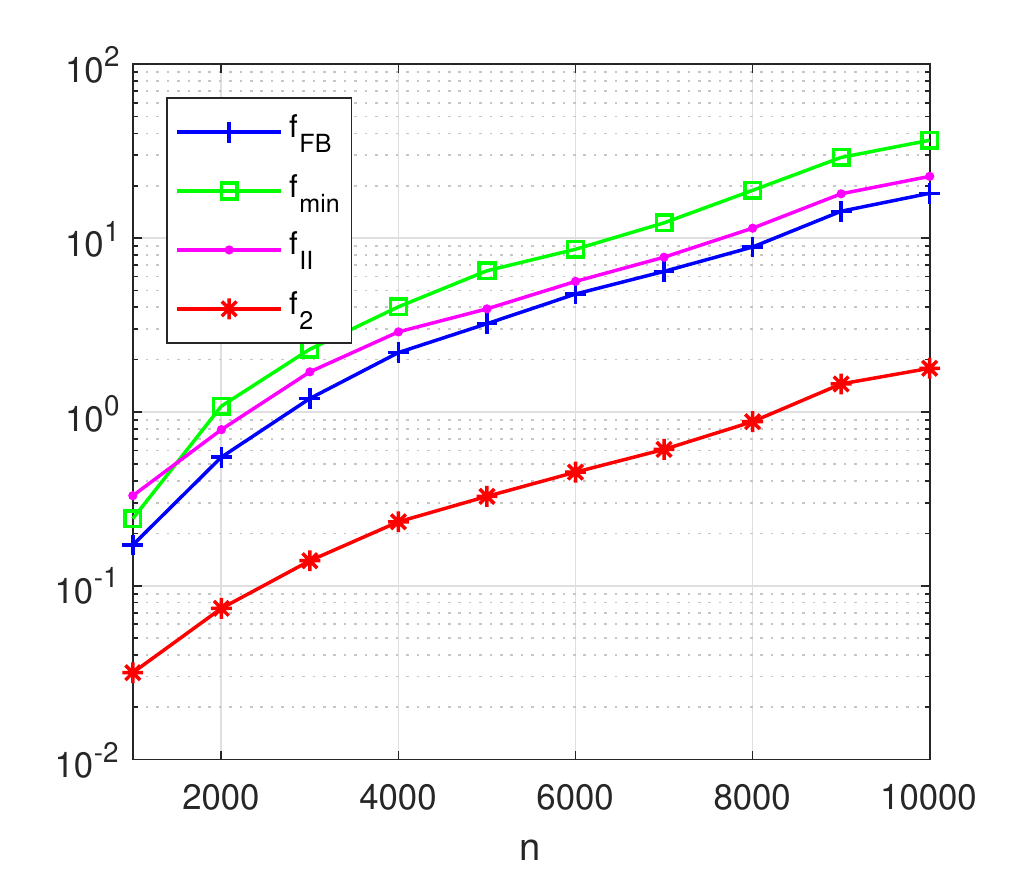}
  \caption{ \Cref{sdp-matrix}: Time}
\end{subfigure} 
\begin{subfigure}{0.32\textwidth}
  \includegraphics[width=.995\linewidth]{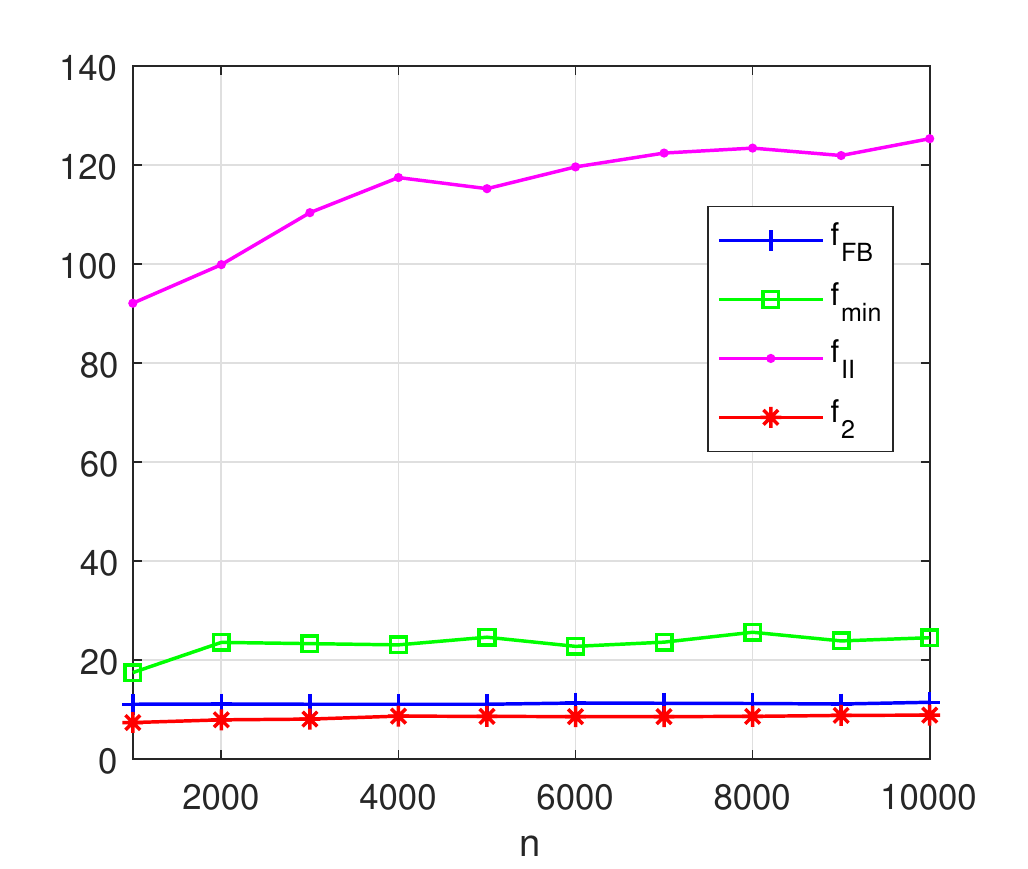}
  \caption{  \Cref{sdp-matrix}: Iteration}
\end{subfigure} \vspace{-5mm}
\caption{Comparison of \NHTP\ solving the sparsity constrained model with four  functions.}
\label{fig:high-dim}
\vspace{-5mm}
\end{figure}

  We now solve the sparsity constrained  model with higher dimensions $n$, and only present  results of  \Cref{z-matrix} and \Cref{sdp-matrix} in   \Cref{fig:high-dim}, since the  results of the rest examples are similar. In terms of accuracy, $f_2$ outperforms the others since it obtains smallest objective function values, with the order of $10^{-17}$ from $f_2$ v.s. $10^{-12}$ from $f_{II}$ in sub-figure (d).  For the computational speed, it can be clearly seen that $f_2$ allows \NHTP\ to run the fastest. By contrast, $f_{FB}$ and $f_{\min}$ run the slowest for  \Cref{z-matrix} and \Cref{sdp-matrix}, respectively. More detailed, {{as  expected}}, $f_{FB}$ and $f_{2}$ for  \Cref{sdp-matrix} in (f) (or $f_{\min}$ and $f_{2}$ for \Cref{z-matrix} in (c)) need similar number of iterations. However, the model with $f_{2}$ makes the method take much shorter CPU time, which means the computational complexity in each step is much lower. Finally, again $f_{II}$ leads to \NHTP\  using more iterations and thus consuming longer total time than that from $f_2$. In summary, among these merit functions, the sparsity constrained model with $f_{2}$ allows \NHTP\ to run the fastest to get the most desirable solutions.

  \section{Conclusion}
A new merit function $f_r$ has been introduced to {{convert the sparse LCP into a sparsity constrained optimization,  enjoying many properties}}, such as being  continuously differentiable for any $r\geq2$,  twice  continuously differentiable for any $r>2$, and convex if the matrix is positive semidefinite. The relationship between the stationary points to the sparsity constrained optimization and solutions to the sparse LCP has been well revealed.  Numerical experiments demonstrated that the adopted  method \NHTP\ has excellent performance to solve the sparsity constrained optimization. {{Most importantly,  comparing the merit functions constructed from other existing famous NCP functions, the optimization with our merit function $f_r$ enables \NHTP\ to possess the lowest computational complexity, fastest convergent speed and most desirable accuracy.  As a result, through converting  the sparse LCP into the sparsity constrained optimization with the help of $f_r$, it can be effectively  solved by \NHTP. In addition, we feel that the new proposed NCP function $\phi_r$ might be able to deal with the sparse nonlinear complementarity problem. We will explore more on this topic in future.}}

\appendix
\section{Proof of theorems in \Cref{sec:pre} - \Cref{sec:nhtp}} 
\subsection{Proof of \Cref{P-matrix}}
The result 1) is taken from \cite[Theorem 3.3.4]{CPS92}. We prove the second claim. If $A$ is a P$_s$-matrix,  then for each nonzero $x\in\R^n$ with $T:=\supp(x)$ and $|T|=\|x\|_0\leq s$, $A_{TT}$ is a $P$ matrix by the definition of  P$_s$-matrix. This implies there is an $i\in T$ such that
$x_i(Ax)_i=x_i(A_{TT}x_T)_i>0.$
Conversely, if for each nonzero $x\in\R^n$ with $T=\supp(x)$ and $|T|\leq s$,  {{then there is an}} $i$ such that $x_i(Ax)_i>0.$ Clearly, such $i\in T$. Since $(Ax)_T=A_{TT}x_T$, this statement is equivalent to that for each given $T$ with $|T|\leq s$, for each nonzero $z\in\R^{|T|}$, {{there is a}} $j$ such that $z_j(A_{TT}z)_j>0$. Therefore, $A_{TT}$ is a P-matrix. Moreover, $T$ can be any subset of $\N$ with $|T|\leq s$, so any $A_{TT}$ is a P-matrix, which means $A$ is a $P_s$ matrix.

\subsection{Proof of \Cref{pro-LCP}} 
1) It follows from  \Cref{pro-phi} that $\phi_r(a,b)$ is continuously differentiable. This together with  $x_i=\langle e_i, x\rangle$ and $y_i=M_ix+q_i$ both being continuously differentiable leads to $\phi_r(x_i,y_i)$ {{being also}} continuously differentiable. Then the $\nabla f_r(x)$ is derived by the addition and chain rules, namely,
\begin{eqnarray*}
\nabla f_r(x)&=&\sum  \left[\partial_1\phi_r(x_i,y_i)\nabla x_i+\partial_2\phi_r(x_i,y_i)\nabla y_i\right]\\
&=&\sum   \left[\left((x_i)_+^{r-1} (y_i)_+^r-|(x_i)_-|^{r-1}\right)e_i+\left((x_i)_+^{r} (y_i)_+^{r-1}-|(y_i)_-|^{r-1}\right) M_i^\top \right]\\
&=&  x_+^{r-1}\circ y_+^r - |x_-|^{r-1} + M^\top \left[x_+^{r}\circ y_+^{r-1} - |y_-|^{r-1}\right].
\end{eqnarray*}
2) For $r>2$,  $\nabla f_r(x)$ is continuously differentiable because all involved functions in $\nabla f_r(x)$ are continuously differentiable. We omitted the detailed calculations here since the addition and chain rules enable us to derive $\nabla^2 f_r(x)$ directly.

3) When $r=2$, it follows
$$\nabla f_2(x)=  x_+ \circ y_+^2 - |x_-|+ M^\top\left[x_+^2\circ y_+ - |y_-| \right]
= \underset{=:g(x)}{\underbrace{x_+ \circ y_+^2 +x_-}}+ \underset{=:h(x)}{\underbrace{M^\top \left[x_+^2\circ y_+ +y_- \right]}}.$$
Then from \cite[Proposition 1.12]{C75} or  \cite[Example 2.6]{HSN84}, we have
$$\partial^2 f_2(x)= \partial (\nabla f_2(x)) \subseteq  \partial  g(x) +  \partial h(x).$$
Therefore, the next step is to calculate  $\partial  g(x)$ and  $\partial h(x)$. For each $g_i(x)$, we have
\begin{eqnarray*}g_i(x)=\left(x_+ \circ y_+^2 +x_-\right)_i=\left\{\begin{array}{ccc}
                                             x_i (y_i)_+^2, & x_i>0, \\
                                             x_i, & x_i\leq0.
                                           \end{array}\right.
\end{eqnarray*}
It is easy to obtain that the generalized Jacobian of $g_i(x)$ by
\begin{eqnarray*}\partial g_i(x)=\left\{\begin{array}{rr}
                                             \left \{ 2x_i (y_i)_+ M^\top_i+(y_i)_+^2 e_i  \right\} , & x_i>0, \\
                                             {\rm co}\left\{ e_i,(y_i)_+^2 e_i\right\}, & x_i=0,  \\
                                             \left\{ e_i\right\} , & x_i<0,
                                           \end{array}\right.
\end{eqnarray*}
which implies that $$\partial g(x)= \{ ~2\Diag(x_+ \circ y_+)M+\Diag(\xi):~\xi\in\Omega_\xi(x) ~\},$$ where $\Omega_\xi(x)$ is given by (\ref{xilcp}). {{Similar reasoning}} also allows us to derive $$\partial h(x)= \{~2 M^\top\Diag(x_+ \circ y_+)+ M^\top \Diag(\zeta) M :~\zeta\in\Omega_\zeta(x)~\},$$ where $\Omega_\zeta(x)$ is given by (\ref{zetalcp}). Those prove the claim.

4) For any $r\geq2$, it follows from (\ref{hess-phi-lcp}) and (\ref{hess-phi2-lcp}) that $\nabla^2f_r(x)$ and any element in $\partial^2f_2(x)$ are positive semidefinite  if $M\succeq 0$ and thus $f_r(x)$ is convex.

\subsection{Proof of \Cref{the-sta}}
1) If $M$ is positive semidefinite and \fea\ is nonempty, it follows from \cite[Theorem 3.1.2]{CPS92} that \sol\ is nonempty. Then $ \sol ={\rm argmin}_{x}  \ f_r(x)$ by (\ref{sol-sco}). Again $M$ being positive semidefinite results in the convexity of $f_r$ {{from \Cref{pro-LCP} 4)}}, which means a point $x^*$ is a solution to ${\min}_{x}  \ f_r(x)$ if and only if $\nabla f_r(x^*)=0$, namely, a stationary point.
 
2)  {{If $M$ is a P-matrix,  we can conclude from \cite[Theorem 5.1, Lemma 5.2]{KYF97} that a point is a solution to (\ref{CP}) if and only if it is a stationary point}}. Thus we have {\sol}$=\mathcal{G}_f$. Then by \cite[Theorem 3.3.7]{CPS92} or \cite[Theorem 1.4]{R07}, (\ref{CP}) has a unique solution for all   $q\in\R^n$ if and only if $M$ is a P-matrix.

 \subsection{Proof of \Cref{exist-s-quad}} The problem (\ref{s-quad}) is equivalent to 
\begin{eqnarray}\label{s-quad-eq}
\underset{T\subseteq\N,|T|\leq s}{{\rm min}}\left\{ \underset{x}{\rm min}~ \langle x, Mx+q \rangle,~~{\rm s.t.} ~ x_T\geq0, x_{T^c}=0,  Mx+q \geq0\right\}.
\end{eqnarray}
Since \feas\ is nonempty, there are some $T$ with $T\subseteq\N,|T|\leq s$ such that the inner program of (\ref{s-quad-eq}) is feasible. This together {{with the Frank-Wolfe theorem \cite{frank1956algorithm} implies that the inner program admits}} an optimal solution $x(T)$ because it is a quadratic program being bounded from below over the feasible region. Clearly, the optimal function value denoted as $\gamma_T$ is unique.  As the choices of $T$ are {{finitely many, e.g., $T\in\{T_1,\cdots, T_N\}$, there are finitely many}} $\gamma_T$. To derive the optimal solution of (\ref{s-quad}), we can pick one $T_i$ from $\{T_1,\cdots, T_N\}$ such that the objective function value  $\gamma_{T_i}$ is the smallest. Then $x(T_i)$ is  an optimal solution of (\ref{s-quad}), namely, $Q_s(M,q)$ is nonempty. \qed
 
 \subsection{Proof of \Cref{existence-main}}\label{proof-exist}
1) Since $M$ is {{ symmetric, $M_T$  having full column rank}} means that $\{M_i^\top:i\in T\}$ are linearly independent. Then it follows from this fact and \cite[Corollary 2.8,  Theorem 3.6]{PXF17}, a global optimal solution $x$ with $\|x\|_0= s$ satisfies the following first order optimality conditions, {{for some $u\in \R^n$,}}
\begin{eqnarray}\label{first-order-11}
\left\{
\begin{array}{r}
 M_{T T}x_T+q_T + M_{TT} x_T-  M_{T\Gamma } u_\Gamma  =0,\\
x_T>0,~x_{T^c}=0,~ u_\Gamma\geq0,~ u_{\Gamma^c}=0,\\
 M_{\Gamma T}x_T+q_\Gamma = 0,~ M_{\Gamma^c T}x_T+q_{\Gamma^c} > 0.
\end{array}
\right.
  \end{eqnarray}
{{where $T$ and  $\Gamma$ are defined in (\ref{TG}).}} We now prove that $T\subseteq \Gamma$. In fact, if there is an $j\in T$ but $j\notin \Gamma$, we have $M_{j T}x_T+q_{j} > 0$ from the last inequality in (\ref{first-order-11}), which derives that $M_{j T}x_T > -q_{j}\geq 0$ by assumption $q_T\leq 0.$ Now consider the first equation in (\ref{first-order-11}),
$$0=M_{j T}x_T+q_j + M_{jT} x_T-  M_{j \Gamma } u_\Gamma
>-  M_{j \Gamma } u_\Gamma  \geq 0$$
due to {{$M$ being a Z-matrix}}, $j\notin \Gamma$ and $u_\Gamma\geq0$. Clearly, this is a contradiction. Therefore, we have $T\subseteq \Gamma$, namely,  $M_{T T}x_T+q_T = 0$, which gives rise to $\langle x, Mx+q\rangle= \langle x_T, M_{T T}x_T+q_T \rangle=0$. Thus $x\in\sol$, showing $x\in\sol\cap S$.

2) Since $M$ is {{symmetric, $M_{T\Gamma}$  having full column rank means that $\{M_{iT}^\top:i\in T\}$ are linearly independent. From  this and \cite[Corollary 2.8, Theorem 3.6]{PXF17}, a global optimal solution $x$ with $\|x\|_0<s$ satisfies the following first order optimality conditions, for some $u,v\in \R^n$,
\begin{eqnarray}\label{first-order-2}
\left\{
\begin{array}{r}
Mx+q+M  x-v-M  u=0,\\
x\geq0,~ v\geq 0,~ \langle x,v \rangle=0,\\
{{u\geq0,~ Mx+q\geq 0,~ \langle u,Mx+q\rangle=0.}}\\
\end{array}
\right.
  \end{eqnarray}
In addition, $v_T=0$ by $x_T>0$ and $\langle x,v \rangle=0$. This and  the above conditions suffice  to \eqref{first-order-11}.  Then the rest of proof is the same as that of proving 1). }}

3) Since $M$ is {{symmetric, $M_{T\Gamma}$  having full column rank means that $\{M_{iT}^\top:i\in T\}$ are linearly independent. By 2), we obtain \eqref{first-order-2} which }} can be rewritten as the conditions that
 are identical to ones presented in \cite[Lemma 3.1.1]{CH97}. Then $M$ being positive semidefinite and \cite[Theorem 3.1.2]{CH97} allow  us to conclude the result.

 \subsection{Proof of \Cref{existence}}
If $|\theta|$=0, then $q\geq0$, which results in $x^*=0$ being a solution to (\ref{sncp}), and thus the conclusion holds immediately. Now consider  $0<|\theta|\leq s$. Clearly, $M_{\theta\theta}$ is a P matrix since $M$ is a P$_s$ matrix. {{This and \Cref{the-sta} 2) allow}} us to conclude that there is a unique solution $u$ satisfying
\begin{eqnarray}\label{u-sol}u\geq 0,~~ M_{\theta\theta}u+q_\theta\geq 0,~~ \langle u, M_{\theta\theta}u+q_\theta\rangle=0.\end{eqnarray}
Since $M\geq0$ and $q_{\theta^c}\geq0$ because of $0<|\theta|\leq s$, we have
$M_{\theta^c\theta}u+q_{\theta^c}\geq 0$. Finally, by letting $x^*_\theta=u$ and $x^*_{\theta^c}=0$, we have $x^*\in(\sol \cap S)$. To see the uniqueness, assume there is another point $z\in(\sol \cap S)$ with  $\supp(z)\subseteq \theta$. Clearly, $z_\theta$ satisfies (\ref{u-sol}). {{However,}}  (\ref{u-sol}) only admits one solution $u$. Therefore, $z_\theta=u=x_\theta$. 
 
 \subsection{Proof of \Cref{pro-bound-s}}
Suppose there is an unbounded subsequence of $\{x^k\}_{k\in K}\subseteq \mathcal{L}_s(f_r,\gamma)$ for
some $\gamma\geq0$, where $K$ is a subset of $\{1,2,\cdots\}$. Let the index set
$J := \{i \in \N:~ \{x_i^k\}~ {\rm is~unbounded}\},$
which is nonempty due to $\{x^k\}_{k\in K}$ {{being unbounded.}} Now define a bounded sequence $\{z^k\}$  by
$$z_i^k=\left\{
\begin{array}{cc}
  0, & i\in J, \\
x_i^k, & i\notin J.
\end{array}
\right.$$
Clearly, we have $z^k\in S$ and $x^k-z^k\in S$ due to $x^k\in S$. Now since $M$ is a P$_s$ matrix (see  \Cref{P-matrix}), then there exists a $\tau>0$ such that $\max_j(z_jM_jz)\geq \tau\|z\|^2$ for each nonzero $z\in S$. In fact, if for any $\tau>0$, there is a nonzero $z\in S$ such that $\max_j(z_jM_jz)< \tau\|z\|^2$, {{then we have}} $z^\top Mz=\sum_j z_jM_jz < n \tau\|z\|^2$, which leads to
$$\tau>\frac{z^\top Mz}{n\|z\|^2}=\frac{z_T^\top M_{TT}z_T}{n\|z_T\|^2}\geq \frac{\sigma_{\min}(M_{TT})}{n}>0,$$
where $\sigma_{\min}(M_{TT})$ is positive due to $M_{TT}$ being a P matrix from $M$ {{being a P$_s$ matrix}}, which is a contradiction if $\tau$ is sufficiently small. So, {{the above assertion indicates}}
\begin{eqnarray*}
\tau\sum_{i\in J} (x_i^k)^2&=&\tau\|x^k-z^k\|^2
\leq \max_j ~( x^k_j-z^k_j)M_j(x^k -z^k)\\
&=& \max_{j\in J} ~(M_jx^k -M_jz^k)( x^k_j-z^k_j)= (M_{j_0}x^k -M_{j_0}z^k)x^k_{j_0}\\
&\leq& (|M_{j_0}x^k|+|M_{j_0}z^k|)|x^k_{j_0}|,
\end{eqnarray*}
where the first inequality {{comes from}} $x^k-z^k\in S$ and $j_0$ is one of the indices for which the max is attained. This inequality divided by $|x^k_{j_0}|$ on both sides derives that
$$\tau|x^k_{j_0}|\leq\tau|x^k_{j_0}|+\tau\sum_{i(\neq j_0)\in J} (x_i^k)^2/|x^k_{j_0}|\leq |M_{j_0}x^k|+|M_{j_0}z^k|.$$
Since $\{z^k\}$ is bounded and $Mx+q$ is continuous, $|M_{j_0}z^k|$ is bounded. Because of this, the above inequalities suffice to $|M_{j_0}x^k|\rightarrow \infty$ as $k(\in K)\rightarrow \infty$. Thus, $|x^k_{j_0}|$ and $|M_{j_0}x^k|$ both tend to infinity, {{leading to $f_r(x^k)\rightarrow \infty$.  Clearly, this  contradicts the definition of the level set that  $f(x^k)\leq\gamma$. }}

Moreover, $\mathcal{O}_s:={\rm argmin}_{x\in S} f_r(x)\subseteq \mathcal{L}_s(f_r,f_r(0))$ is bounded as the level set is bounded. If $(\sol \cap S)=\emptyset$, then the conclusion holds readily. If  $(\sol \cap S)$ is nonempty, then for any $x^*\in (\sol \cap S)$  it follows $f_r(x^*)=0$, which means $x^*\in \mathcal{O}_s$ due to  $f_r(x)\geq 0$. Namely,  $(\sol \cap S)\subseteq\mathcal{O}_s$.
 
 \subsection{Proof of \Cref{the-sta-s}}
It follows from (\ref{grd-phi-lcp}) that
\begin{eqnarray}\label{nabla-f}\nabla f_r(x)=   x_+^{r-1}\circ y_+^r - |x_-|^{r-1}  + M^\top \Big[ x_+^{r}\circ y_+^{r-1} - |y_-|^{r-1} \Big],\end{eqnarray}
where $y:=Mx+q.$ If $x$ is a solution to (\ref{sncp}), namely, $x\geq0, y\geq0,\langle x,y\rangle=0$ and $\|x\|_0\leq s$, then $x$ is {{a stationary point}} due to $\nabla f_r(x)=0$ satisfying (\ref{necc-opt-cod-1}). We now  prove the second part. For any $x$ with $T=\supp(x)$ {{such that (\ref{necc-opt-cod-1}) holds}}, besides $T_+$ and  $\Gamma_+$, let
\begin{eqnarray}\label{sets}
\begin{array}{rll}
T_-&:=&\{i\in\N:~x_i<0\}, ~~\Gamma_-~:=\{i\in\N:~y_i<0\}, \\
\alpha&:=&T_+\cap\Gamma_+~~=\{i\in\N:~x_i>0, y_i>0\}, \\
\beta&:=&T_+\setminus \alpha ~~~~=\{i\in\N:~x_i>0, y_i\leq 0\}.
\end{array}
\end{eqnarray}
Clearly, $T=T_-\cup\alpha\cup \beta$. From (\ref{necc-opt-cod-1}), {{$x$ is a stationary point, then $\nabla_{T} f_r(x)=0$. Based on the above notation}},  (\ref{nabla-f}) allows us to  write $\nabla_{\alpha} f_r(x)$ as
\begin{eqnarray}
 0=\nabla_{\alpha} f_r(x)&=&  (x_{\alpha})_+^{r-1}\circ (y_{\alpha})_+^r - |(x_{\alpha})_-|^{r-1}  + M_{\alpha}^\top  [ x_+^{r}\circ y_+^{r-1} - |y_-|^{r-1}  ],\nonumber \\
&=&x_{\alpha}^{r-1}\circ y_{\alpha}^r + M_{\alpha \alpha}^\top (x_{\alpha}^r\circ y_{\alpha}^{r-1})- M_{\Gamma_-\alpha}^\top|y_{\Gamma_-}|^{r-1}\nonumber\\
&\geq& x_{\alpha}^{r-1}\circ y_{\alpha}^r + M_{\alpha \alpha}^\top (x_{\alpha}^r\circ y_{\alpha}^{r-1})\nonumber\\
\label{g-a} &= &  \left(\Diag(y_\alpha)+ M_{\alpha \alpha}^\top\Diag(x_\alpha)\right) (x_{\alpha}^{r-1}\circ y_{\alpha}^{r-1})=:A(x_{\alpha}^{r-1}\circ y_{\alpha}^{r-1}), 
 \end{eqnarray}
 where {{the inequality}} holds due to $\Gamma_-\cap\alpha\neq\emptyset$ and $- M_{\Gamma_-\alpha}^\top|y_{\Gamma_-}|^{r-1}\geq 0$ by $M$ being a Z matrix.
If $\alpha\neq\emptyset$, then $x_{\alpha}>0,~ y_{\alpha} >0$ and $A\succ0$ due to $M_{TT}\succeq0$ and $\alpha\subseteq T_+$.  Multiplying  both sides of (\ref{g-a}) by $\nu:=(x_{\alpha}^{r-1}\circ y_{\alpha}^{r-1})^\top$ derives
$0\geq \nu^\top A \nu>0,$  which clearly is a contradiction. {{Thus $\alpha = \emptyset$, giving}} rise to   $x_{+}\circ y_{+} =0 $ and $T_+=\beta$. Now, $T=T_-\cup T_+$ and $\nabla_{T} f_r(x)=0$ leading to
\begin{eqnarray}\label{g-a-b-t}
 &&~~~~0= 
\left[
\begin{array}{r}
M_{\Gamma_-\beta}^\top|y_{\Gamma_-}|^{r-1}\\
|x_{T_-}|^{r-1}+ M_{\Gamma_-T_-}^\top|y_{\Gamma_-}|^{r-1}
\end{array}
\right]= \left[
\begin{array}{rr}
0&M_{\Gamma_-T_+}^\top\\
I&  M_{\Gamma_-T_-}^\top
\end{array}
\right]  \left[
\begin{array}{r}
|x_{T_-}|^{r-1}\\
|y_{\Gamma_-}|^{r-1}
\end{array}
\right] =:Bz.
 \end{eqnarray}
Clearly, $z>0$ from {{the definitions}} of $\Gamma_-$ and $T_-$.   {{Stiemke  Theorem (see \cite[Theorem 13]{P17} or \cite[Theorem 7]{M69})}} states that $Bz=0,z>0$ has no solution if  $B^\top u\geq0, u\neq0$ has a solution. By assumption,  there is a nonzero $v\in\R^{|T_+|}$ such that $M_{\Gamma_+^cT_+}v\geq0$, which indicates  $M_{\Gamma_-T_+}v\geq0$ due to $\Gamma_-\subseteq\Gamma_+^c$. Let $u=[v^\top~0]^\top\neq0$, then  we have
$B^\top u=  [0 ~ (M_{\Gamma_-T_+}v)^\top]^\top \geq0.$
Thus $Bz=0, z>0$ has no solution, which implies that  $z=0$ and hence $\Gamma_-=  T_-=\emptyset$. Those {{together with $\alpha=0$ enable us to obtain $ x\geq 0, y\geq0, x\circ y=0$}}. Finally, it follows from $x\in S$ owing to $x$ satisfying (\ref{necc-opt-cod-1}) that $x\in\sol \cap S$.  

 \subsection{Proof of \Cref{theorem-local-stationary}}
    1) The sufficiency is derived by (\ref{necc-opt-cod}) and (\ref{necc-opt-cod-1}) easily. We now prove the necessity. Since $M$ is positive semidefinite, $f_r$ is a convex function from {{Lemma \ref{pro-LCP} 4).}} As $x^*$ is a stationary point (\ref{necc-opt-cod-1}) with $\|x^*\|_0<s$, $\nabla f_r(x^*)=0$.  Then for any $x\in\R^n$, it holds
 \begin{eqnarray}\label{global}f_r(x)\geq f_r(x^*)+\langle \nabla f_r(x^*), x-x^*  \rangle=  f_r(x^*),  \end{eqnarray}
    which shows the global optimality of $x^*$. If further  \texttt{fea}$(M,q)$ is nonempty, then $\sol$ is nonempty from   \Cref{the-sta} 1). Now replacing $x$ by any $z\in\sol$ in (\ref{global}) yields $0=f_r(z)\geq  f_r(x^*)\geq 0$, which means $x^*\in\sol$ and hence $x^*\in(\sol\cap S)$.

    2) The sufficiency is obvious by (\ref{necc-opt-cod}) and (\ref{necc-opt-cod-1}).   By (\ref{necc-opt-cod-1}), $x^*$ being a stationary point with $\|x^*\|_0=s$ leads to $\nabla_{T_*}  f_r(x^*)=0$.  Then for any $x\in\R_{T_*}$, we have
     \begin{eqnarray*}f_r(x) \geq f_r(x^*)+\langle \nabla f_r(x^*), x-x^*  \rangle=f_r(x^*)+\langle \nabla_{T_*} f_r(x^*), x_{T_*}-x^*_{T_*}  \rangle = f_r(x^*). \end{eqnarray*}
    This proves the local optimality of $x^*$. If  $M_{T_*T_*}$ is nonsingular, then (\ref{zetalcp}) yields
  $$  {{\begin{array}{l}
    \nabla^2_{T_*T_*}f_2(x^*)\succeq M_{T_*T_*}\Diag(\zeta_{T_*})M_{T_*T_*}~~{\rm with}~~~~\zeta_i\in\Xi(y_i,x_i).\end{array}}}$$
Clearly,  $\zeta_{T_*}>0$ due to $x_i\neq0, i\in T_*$ and hence $\nabla^2_{T_*T_*}f_2(x^*)\succ \lambda I$, where $\lambda$ is the smallest eigenvalue of $(M_ {T_*T_*}\Diag(\zeta_{T_*})M_{T_*T_*})$. Then for any  $x\in\R_{T_*}$, it holds
\begin{eqnarray*}f_2(x)&\geq& f_2(x^*)+\langle \nabla f_2(x^*), x-x^*  \rangle+(\lambda/2)\|x-x^*\|^2> f_2(x^*),  \end{eqnarray*}
which shows the global optimality of $x^*$ on $\R_{T_*}$. 
 \subsection{Proof of \Cref{bound-c}} Since $M$ is a  P$_s$ matrix, then $\mathcal{L}_s(f_r,f_r(0))$ is bounded from \Cref{pro-bound-s} and thus $x$ is bounded, which suffices {{to the boundedness of $y:=Mx+q$.}}  By  (\ref{hess-phi-lcp}) we conclude that $\nabla^2 f_r(x)$  is bounded  for any $r>2$. For $r=2$, from (\ref{hess-phi2-lcp}), any point in $\partial^2 f_2(x)$ is bounded since both  $\Omega_\xi(x)$ and $\Omega_\zeta(x)$ are bounded. Namely, $\nabla^2 f_2(x)$  is bounded as well. Therefore, there exists $C<+\infty$  such that $\sigma_{\max}(\nabla^2 f_r(x))< C$ for any $x\in\mathcal{L}_s(f_r,f_r(0))$.
 
    \subsection{Proof of \Cref{theorem-convergence}}
1) Choice of $x^0\in\mathcal{L}_s(f_r,f_r(0))$ indicates that $\nabla^2  f_r(x^0)\preceq C I_n$ by \Cref{bound-c}. This together {{with the reasoning}} to prove Lemma 5 in \cite{ZXQ19}, in which we set $T_{-1}\supseteq\supp(x^0)$ with $|T_{-1}|=s$ and replace $M_{2s}$ by $C$,  derives
\begin{equation}
\langle d^0, \nabla   f_r(x^0)\rangle \leq -\rho \|d^0\|^2-(\eta/2) \|\nabla_{T_{-1}}   f_r(x^0)\|^2,
\end{equation}
where $\rho>0$ is a constant associated with $\mu$ and $ C$. Then {{the same reasoning}} to proof Lemma 7 in \cite{ZXQ19} derive that
\begin{equation}
 f_r(x^1)- f_r(x^0)  \leq -\rho_1 \|d^0\|^2-(\eta_1/2) \|\nabla_{T_{-1}}   f_r(x^0)\|^2\leq 0,
\end{equation}
where $\rho_1>0,\eta_1>0$ are two constants associated with $\mu$ and $ C$. So, $f_r(x^1)\leq f_r(x^0)\leq f_r(0)$, which means $x^1\in\mathcal{L}_s(f_r,f_r(0))$ and because of this, $\nabla^2  f_r(x^1)\preceq C I_n$. In addition, $T_{0}\supseteq\supp(x^1)$ with $|T_{0}|=s$ from \Cref{alg:nhtp}.  By the induction, we can conclude that
\begin{equation}
 f_r(x^{k+1})- f_r(x^k)  \leq -\rho_1 \|d^k\|^2-(\eta_1/2) \|\nabla_{T_{k-1}}   f_r(x^k)\|^2\leq 0,
\end{equation}
for any $k=0,1,2,\ldots.$ This displays  the non-increasing property of $\{f_r(x^k)\}$ and {{derives}}  $f_r(x^k)\leq f_r(x^0)\leq f_r(0)$. Consequently, $x^k\in\mathcal{L}_s(f_r,f_r(0))$ and {{it is bounded.
The  proofs of 2) and 3) are the same}} as those of proving Lemma 7, Theorem 8 and Theorem 9 in \cite{ZXQ19}. We omit them here. 

\section*{Acknowledgments}
We sincerely thank the associate editor and the two referees for their
detailed comments that have helped us to improve the paper.  We also thank Prof. Naihua Xiu of Beijing Jiaotong University who offered us valuable instructions.

\bibliographystyle{siamplain}
\bibliography{references}
\end{document}